\documentclass[12pt,a4paper]{article}
\usepackage{amsmath,amssymb,amsfonts,epic,epsfig,graphics}
\usepackage[T1]{fontenc}
\usepackage[utf8]{inputenc}
\usepackage{epsfig}
\usepackage{faktor}
\usepackage[all]{xy}
\usepackage[german]{babel}
\newcommand{\ko}{{\mathcal O}}

\newcommand{\kh}{{\mathcal H}}
\newcommand{\kf}{{\mathcal F}}
\newcommand{\Z}{\mathbb{Z}}

\renewcommand{\P}{\mathbb{P}}

\newcommand{\C}{{\mathbb C}}
\newcommand{\R}{{\mathbb R}}
\newcommand{\SL}{{\text{SL}}}
\newcommand{\GL}{{\text{GL}}}

\begin{document}
\vspace{1cm}

\begin{center}
\textbf{\large Leben und Werk von Egbert Brieskorn (1936 -- 2013)}


\textbf{\em  Gert-Martin Greuel, Walter Purkert}
\end{center}

\begin{center}
\includegraphics[width=7cm,height=8cm]{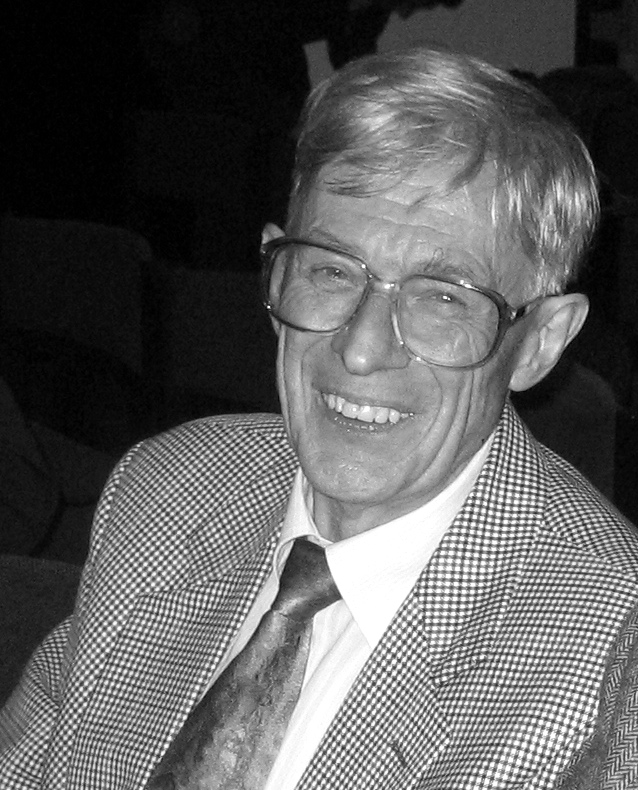}\\
Brieskorn 2007
\end{center}

\bigskip

Egbert Brieskorn verstarb am 11. Juli 2013, wenige Tage nach seinem 77. Geburtstag. Er war eine beeindruckende Persönlichkeit, die auf alle, die ihn näher kannten, sei es innerhalb oder außerhalb der Mathematik, einen nachhaltigen Eindruck hinterlassen hat. Brieskorn war ein großartiger Mathematiker, aber seine Interessen, sein Wissen und seine Aktivitäten reichten weit über die Mathematik hinaus. In dem folgenden Beitrag, der stark von einer viele Jahre andauernden persönlichen Verbundenheit der Autoren mit Brieskorn geprägt ist, versuchen wir, einen tieferen Einblick in das Leben und Werk von Brieskorn zu geben. Dabei beleuchten wir sowohl sein persönliches Engagement für Frieden und Umwelt als auch seine langjährige Erforschung des Lebens und Werkes von Felix Hausdorff und die Herausgabe von Hausdorffs Gesammelten Werken; der Schwerpunkt  des Artikels liegt jedoch auf der Darstellung seines bemerkenswerten und einflussreichen mathematischen Werkes.

Der erste Autor (GMG) hat prägende Teile seiner wissenschaftlichen Laufbahn als Diplomand und Doktorand mit Brieskorn in Göttingen und später als sein Assistent in Bonn verbracht und schildert in den beiden ersten Teilen, zum Teil aus der Erinnerung der persönlichen Zusammenarbeit, Aspekte aus Brieskorns Leben und von seinem politisch–gesellschaftlichem Engagement. Er erläutert zudem im Teil über Brieskorns mathematisches Werk ausführlich die wissenschaftlichen Hauptergebnisse seiner Publikationen. Der zweite Autor (WP) hatte vor allem im Zusammenhang mit dem Hausdorff-Projekt viele Jahre mit Brieskorn zu tun; der entsprechende Abschnitt über das Hausdorff--Projekt wurde von ihm verfasst.

Wir danken Wolfgang und Bettina Ebeling, Helmut Hamm, Thomas Peternell, Anna Pratoussevitch und Wolfgang Soergel für nützliche Hinweise und vor allem  Brieskorns Frau Heidrun Brieskorn für die Überlassung von Material aus Brieskorns Nachlass.

\bigskip

\textbf{Mathematical Subject Classification (MSC2010)}:  01A61, 14B05, 14B07, 14D05, 14F45, 14H20, 14J17, 14J70, 17B22, 32S05, 32S25, 32S40, 32S55, 57R55\\

\textbf{\large Stationen seines Lebens}
 
Brieskorn wurde am 7. Juli 1936 in Rostock als Sohn eines Mühlenbau–Ingenieurs geboren und wuchs zusammen mit seiner Schwester bei seiner Mutter in Siegerland auf. Über seine Jugend und woher seine Begeisterung für die Mathematik rührte ist nur  wenig bekannt. Immerhin wissen wir aus dem Kapitel \glqq Childhood and Education\grqq\ aus dem Film der Simons Foundation über Brieskorn \cite {EB2010}, dass seine Mutter seine kindliche Neugierde unterstützte und dass sein Vater sein technisches Interesse förderte. Er hatte auch einen guten Mathematiklehrer im Gymnasium, der ihm über den Schulstoff hinausgehende mathematische Literatur besorgte, wobei ihn vor allem geometrische Konstruktionen interessierten (und z.B weniger eine Arbeit von Gauß). Auch wenn sein technisches Interesse zunächst noch überwog, war das Interesse an der Mathematik aber schon vor seinem Studium stark ausgeprägt. Bei dem Aufnahmegespräch für das Evangelische Studienwerk Villigst, die Hochbegabten--Förderungseinrichtung der evangelischen Kirche, meinte der Prüfer zu ihm sinngemäß: "`Herr Brieskorn, Ihr Talent und Ihre Begeisterung für die Mathematik sind außergewöhnlich, aber vergessen Sie nicht, dass es neben der Mathematik noch andere Dinge im Leben gibt."' Diese Episode erzählte Egbert Brieskorn dem erstgenannten Autor dieser Zeilen mit leicht ironischem Unterton viel später, als tatsächlich andere Dinge als die Mathematik sein Leben und Wirken bestimmten. Der Leiter des Evangelischen Studienwerks erkannte auch, dass sein ursprünglicher Studienwunsch - Elektrotechnik - nicht seinem Wesen entsprach und er überzeugte ihn, etwas Theoretisches zu studieren. 

Brieskorn begann daher im Oktober 1956 Mathematik und Physik in München zu studieren. Er wechselte zum Sommersemester 1959 nach 5 Semestern auf Anraten von Karl Stein nach Bonn, um dort den Satz von Hirzebruch–Riemann-Roch zu verstehen, den er als "`my first love in mathematics"' \cite{EB2000} bezeichnete. Friedrich Hirzebruch, der selbst erst 1956 nach Bonn gekommen war, beeindruckte den jungen Studenten Brieskorn durch seine freundliche, offene Persönlichkeit und durch seinen klaren Vortragsstil zutiefst. Brieskorn wurde Hirzebruchs Student und promovierte bei ihm 1963 mit der Dissertation  "`Differentialtopologische und analytische Klassifizierung gewisser algebraischer Mannigfaltigkeiten"'. Hirzebruch bezeichnete Brieskorn später als seinen talentiertesten Schüler und Brieskorn hat seinen Lehrer Hirzebruch sein Leben lang als Mathematiker und Mensch hoch verehrt. 
1968 habilitierte sich Brieskorn in Bonn mit der Habilitationsschrift "`Singularitäten komplexer Räume"' und wurde 1969 als ordentlicher Professor nach Göttingen berufen, wo er bis 1973 blieb. Wegen seiner Frau Heidrun, die er 1973 heiratete und die eine Stelle als Bratschistin am Kölner Rundfunk-Sinfonie-Orchester 
(heute WDR Sinfonieorchester Köln) erhielt, wechselte er 1973 nach Bonn, zunächst an den Sonderforschungsbereich Theoretische Mathematik und ab 1975 auf eine Stelle als ordentlicher Professor, wo er bis zu seiner Emeritierung 2001 wirkte.

\begin{center}
\includegraphics[width=12cm,height=7.5cm]{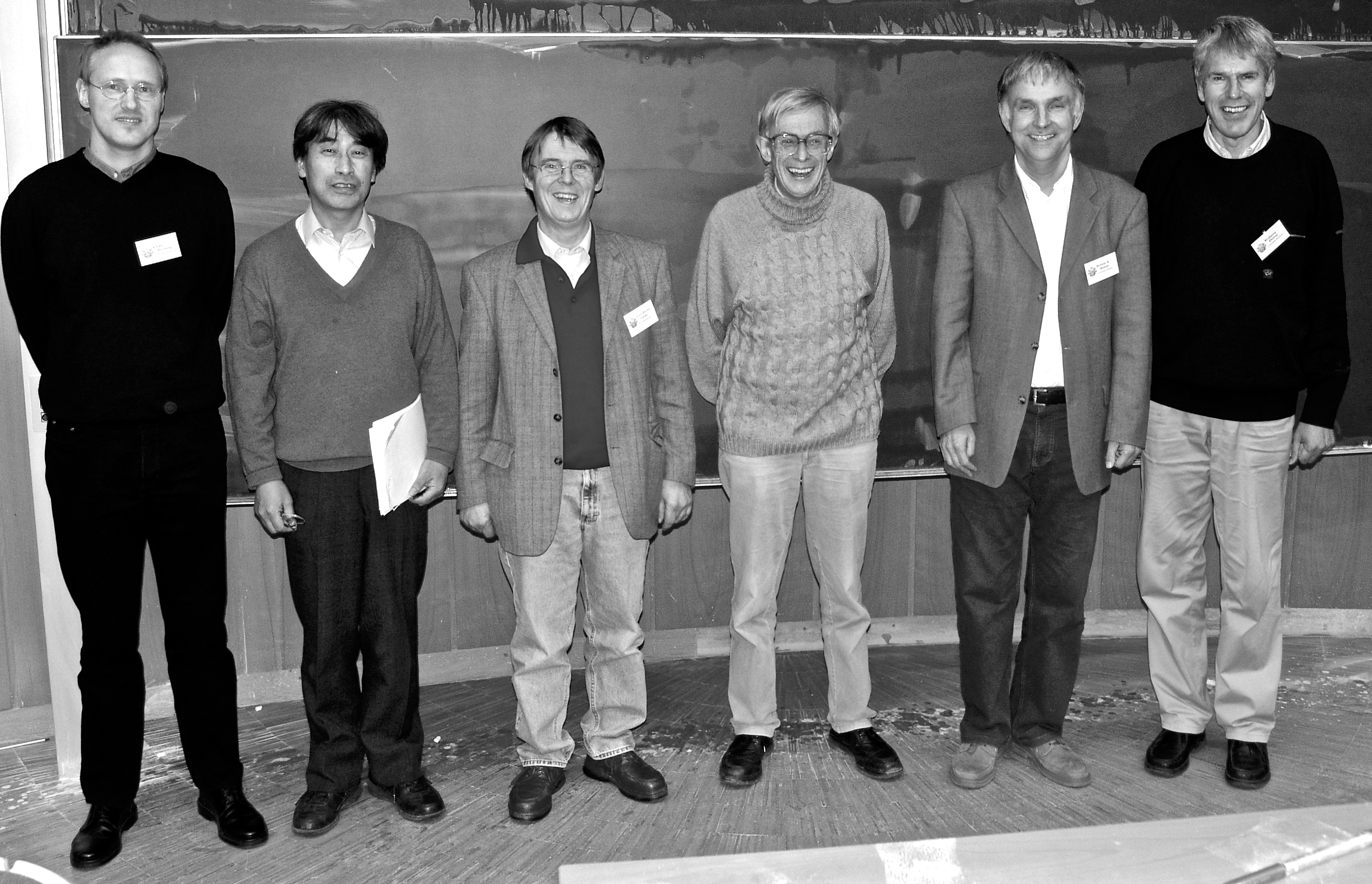}\\

Brieskorn (3. v.r.) und seine Studenten (v.l.) Claus Hertling, Kyoji Saito, Gert-Martin Greuel, Helmut Hamm, Wolfgang Ebeling, 2004
\end{center}

\bigskip

\textbf{\large Mathematik und politisch–gesellschaftliches Engagement}

Obwohl Brieskorn sich schon immer sozial engagiert hatte, z.B. durch Arbeit in einem Stahlwalzwerk der Dortmunder Hörder Hüttenunion im Rahmen eines Werksemesters des Evangelischen Studienwerks, war für ihn während des Studiums und auch noch viele Jahre als Professor die Mathematik das Wichtigste in seinem Leben, wie er selbst schreibt (Lebenslauf Werksemester). Er war vollkommen begeistert von der faszinierenden Schönheit und Klarheit und dem hohen Anspruch der Mathematik. 
Es zeigte sich aber auch in kleinen Dingen, dass er überall mathematische Phänomene sah und von diesen fasziniert war.  Diese Begeisterung für die Mathematik und sein Forschungsgebiet, die Singularitätentheorie, übertrug sich auf die Studierenden. Die folgenden Abschnitte zeigen die Entwicklung des Verhältnisses von Mathematik und politisch–gesellschaftlichem Engagement bei Brieskorn aus der Sicht und Erinnerung des erstgenannten Autors.
Er hatte mich im Spätsommer 1969 zum Tee in seine Göttinger Wohnung eingeladen, um mir das Thema für meine Diplomarbeit zu erklären. Bevor er anfing, entdeckte er eine durch Lichtbrechung entstandene Kaustik in seiner Teetasse, die er als "`einfache Singularität"' klassifizierte, und anschließend interpretierte er die spiralförmige Schokoladenspur in den Keksen als dynamisches System. Er erklärte mir voller Begeisterung exotische Sphären und wie sie sich durch reell–analytische Gleichungen als Umgebungsrand gewisser isolierter Hyperflächensingularitäten beschreiben lassen. Ich war angesteckt, und als er mir dann vorschlug, die Ergebnisse seiner noch unveröffentlichten Arbeit "`Die Monodromie der isolierten Singularitäten von Hyperflächen"' \cite{EB1970a} auf isolierte Singularitäten von vollständigen Durchschnitten zu verallgemeinern, sagte ich sofort zu.

Dass ich überhaupt Brieskorns Student wurde, geschah auf Empfehlung von Hans Grauert, bei dem ich 1966/67 in Göttingen die Anfängervorlesung "`Differential- und Integralrechnung” gehört hatte. Nach meiner Rückkehr von einem einjährigen Studienaufenthalt an der ETH Zürich erhielt ich einen Anruf von Brieskorn, ob ich Interesse hätte, bei ihm eine Diplomarbeit zu schreiben. Grauert hätte mich empfohlen, da wir doch in der gleichen "`Verbindung"' seien (gemeint war das Evangelische Studienwerk). So kam es, dass ich Brieskorns erster Diplomand wurde. 
Brieskorn war im Sommersemester 1969 nach Göttingen gekommen, zunächst zur Vertretung der Stelle und ab Juli 1969 als ordentlicher Professor. Im WS 1969/70 war er zu einem Forschungsaufenthalt an das IHES in Bures-sur-Yvette beurlaubt. Ich hatte seine erste Vorlesung im Sommer 1969 in Göttingen über 2–dimensionale Schemata verpasst und über Winter lernte ich im Lesesaal des Mathematischen Instituts an der Bunsenstraße mit Hilfe von Godements "`Topologie Algébrique et Théorie des Faisceaux"' Garbentheorie und Hyperkohomologie. Die Atmosphäre war dazu hervorragend geeignet. Bis auf gelegentliches Getuschel und das unverwechselbare Geräusch, wenn Carl Ludwig Siegel schwer atmend in den hinteren Teil des Lesesaals zu den alten Werken verschwand, war man ungestört und konnte in dem angrenzenden Diskussionsraum mit Kommilitonen diskutieren. 

Nach seiner Rückkehr aus Frankreich hielt Brieskorn folgende Vorlesungen: "`Differentialtopologie"', dann die Anfängervorlesung "`Infinitesimalrechnung I und II"', anschließend "`Analysis auf Mannigfaltigkeiten"', "`Quantitative Theorie dynamischer Systeme"', "`Algebraische Topologie II"' und vor seinem Wechsel nach Bonn im Sommersemester 1973 "`Einfache Singularitäten"'. Brieskorns Vorlesungen unterschieden sich grundsätzlich von denen von Grauert. Während Grauert nur das erzählte, was er an die Tafel schrieb und bewies, stellte Brieskorn oft größere Zusammenhänge dar und erwähnte Querverbindungen, ohne diese zu beweisen. Offensichtlich hat beides seine Vorteile und die Studierenden in Göttingen schätzten sowohl Brieskorns als auch Grauerts Vorlesungen sehr. Brieskorn betrieb oft einen außergewöhnlichen Aufwand bei der Vorbereitung seiner Vorlesungen, um mathematische wie historische Hintergründe oder Seitenäste des vermittelten Stoffes aufzuzeigen, was man auch an seinen Lehrbüchern deutlich sehen kann. Man kann sagen, dass sein Streben nach Perfektion, das sich im Laufe der Jahre sogar noch steigerte, charakteristisch für ihn war. Perfektion erwartete er auch von seinen Studenten, die er umgekehrt sehr förderte, sowohl mathematisch, u.a. mit wöchentlichen Arbeitstreffen oder mit Empfehlungen, aber auch im persönlichen Bereich. So ließ er z.B. seinen ersten Doktoranden Helmut Hamm während seines Frankreich-Aufenthalts unentgeltlich in seinem Appartement wohnen.

Die Göttinger Zeit von 1969 bis 1973 bewirkte einen bemerkenswerten Wandel in Brieskorns Ansichten und Einstellungen. Neben der Mathematik gewannen auch politische und gesellschaftliche Fragen Bedeutung für ihn. Es war die Zeit kurz nach den heftigen Studierendenprotesten der 68er–Generation, bei der auch die Göttinger Studierenden, in erster Linie die Theologen, danach gleich die Mathematiker (z.B. in einer linken "`Basisgruppe Mathematik"')  sehr aktiv waren. Brieskorn stand einigen Forderungen der Studenten sehr positiv, anderen eher kritisch gegenüber. Er war gegen wissenschaftlich unbegründete hierarchische Strukturen, engagierte sich für stärkere Mitbestimmung des Mittelbaus und der Studentenschaft, und er sympathisierte mit Ideen der Reformuniversitäten in Bremen und Osnabrück. Die wissenschaftliche Qualität des Studiums stand allerdings immer an erster Stelle.
 
Noch wichtiger als studentische Reformvorstellungen waren für Brieskorn jedoch sein Engagement gegen den Vietnamkrieg und damit verbunden sein Eintreten für die Unterstützung unterdrückter Völker. Aus \cite{EB2010} wissen wir, dass er schon während seines Aufenthaltes am Massachusetts Institute of Technology (MIT) zusammen mit Michael Artin und anderen Kollegen vom MIT in New York an einer großen Demonstration gegen den Vietnamkrieg teilnahm, auf der auch Martin Luther King sprach.
So engagierte er sich in Göttingen von Beginn an im Komitee für die wissenschaftliche Zusammenarbeit mit Kuba (KoWiZuKu), das 1970 gegründet wurde und dessen erster Generalsekretär der Mathematiker Klaus Krickeberg war. Dass er auch bei diesen Aktivitäten stets sehr genau war, zeigt eine Episode, als er gemeinsam mit mir und meiner Frau in Göttingen Plakate zur Ankündigung eines Vortrags der kubanischen Gesundheitsministerin in der Universität klebten. Brieskorn achtete penibel darauf, dass alles korrekt war und z.B. keine Plakate an Verteilerhäuschen geklebt wurden, da das zu Hitzeproblemen führen könnte.
Später in Bonn engagierte sich Brieskorn von Beginn an in der Friedensbewegung, die sich aus Protest gegen den "`NATO-Doppelbeschluss"' und damit gegen die Stationierung von Mittelstreckenraketen in der damaligen BRD neu formierte. Brieskorn gehörte zu den Erstunterzeichnern des \glqq Mainzer Appell\grqq, der Abschlusserklärung des Kongresses "`Verantwortung für den Frieden -- Naturwissenschaftler warnen vor neuer Atomrüstung"' im Juni 1983 in Mainz, an der auch der Mathematiker Stephen Smale teilnahm (siehe \cite{HD1983}). In der  "`Naturwissenschaftler-Initiative Verantwortung für den Frieden"' (heute: "`NaturwissenschaftlerInnen-Initiative für Frieden und Zukunftsfähigkeit"') organisierten sich vor allem Physiker, aber auch viele Mathematiker, wie z.B. Brieskorn.

Brieskorns politische Ansichten in der Bonner Zeit sind sicher links von sozialdemokratischen Ideen einzuordnen, sie wandelten sich aber im Laufe der Jahre zu radikal ökologischen Überzeugungen, die er zusammen mit seiner Frau Heidrun lebte. Beide wohnten seit 1982 in einem allein am Waldrand gelegenen Haus in Eitorf an der Sieg, wo sich beide zusammen in einer kleinen Gruppe intensiv dem Naturschutz, genauer dem Artenschutz widmeten. Die Arten, um die es ging, waren zunächst einheimische Fledermäuse, für die sie Winterquartiere in alten Bergwerkstollen kontrollierten, sicherten und neue Winterquartiere bauten. Zur Unterscheidung der verschiedenen Arten machte Heidrun Brieskorn viele Tonaufnahmen, die Egbert Brieskorn mit selbstgeschriebenen Programmen einer Fourier--Analyse unterzog. Noch mehr Arbeitszeit und Mühe verwendeten beide auf den Schutz und die Erhaltung der Lebensbedingungen einer sehr seltenen Schmetterlingsart, des Wiesenknopf--Ameisenbläulings oder lateinisch Maculinea. Dies ging so weit, dass Brieskorn und seine Frau dafür sorgten, dass die Gemeinde ihre Bebauungspläne änderte und sie selbst Grünland aufkauften, das sie selbst von Hand pflegten, um den Lebensraum dieser Falter zu erhalten. Sie gründeten die Maculinea-Stiftung NRW, damit die Arbeit zur Erhaltung der Schmetterlingsarten dauerhaft weitergehen kann. Für ihr Engagement erhielten beide 2013 das \textit{Verdienstkreuz am Bande der Bundesrepublik Deutschland}.

Neben der ehrenamtlichen Arbeit im Naturschutz waren die letzten 20 Jahre von Brieskorns Leben durch die Mitarbeit beim Editionsprojekt "`Felix Hausdorff – Gesammelte Werke"' der Nordrhein-Westfälischen Akademie der Wissenschaften und der Künste bestimmt. Brieskorn selbst schreibt im Juni 2011 kurz vor seinem 75. Geburtstag und wissend, dass er vielleicht nicht mehr lange leben würde, in einem Brief "`An meine lieben ehemaligen Studenten und Schüler"' wie es dazu kam und wie sehr ihn die Biographie von Hausdorff beschäftigt hat: 
"`Eine [der Aufgaben] davon entwickelte sich daraus, dass das Mathematische Institut in Bonn im Januar 1992 den 50.Todestag von Felix Hausdorff feiern wollte. Da kein Kollege einen Vortrag über sein Leben halten wollte, habe ich damals diese Aufgabe übernommen, nicht ahnend, was ich mir damit aufgeladen hatte. Ich habe 20 Jahre damit verbracht, in Archiven und Quellen der verschiedensten Art nach Spuren dieses überaus merkwürdigen Menschen und Mathematikers zu suchen. Ich habe dabei sehr viel gelernt, aber sehr viel Zeit geopfert. Ich habe drei Anläufe unternommen, seine Biographie zu schreiben, und mit der dritten Fassung bin ich, glaube ich, auf dem richtigen Wege. Ungefähr 262 Seiten sind bis jetzt geschrieben, aber noch nicht einmal die Hälfte des Lebensweges dargestellt. Diese Biographie soll im ersten Band einer insgesamt neunbändigen Ausgabe der Werke Felix Hausdorffs erscheinen. Sechs Bände sind gedruckt, mindestens einer, vielleicht aber auch zwei der fehlenden Bände werden in diesem Jahr in Druck gehen. Das Schlimme ist, dass der erste Band, in dem meine Biographie Hausdorffs erscheinen soll, unmöglich rechtzeitig fertig werden kann."' Über die umfangreichen Arbeiten Brieskorns und seine weit über das Übliche hinausgehenden Recherchen zu diesem Projekt wird im folgenden Abschnitt berichtet.
\bigskip

\newpage

\textbf{\large Das Hausdorff-Projekt}

     Die Arbeit, die Egbert Brieskorn in seinem Brief erwähnt, war für mehr als 20 Jahre eine Herzensangelegenheit für ihn, nämlich die Erforschung des Lebens und  Werkes von Felix Hausdorff (1868-1942) und die Herausgabe von Hausdorffs Gesammelten Werken. 

    Felix Hausdorff hat die allgemeine Topologie als eigenständiges Gebiet der Mathematik begründet und darüber hinaus grundlegende Beiträge zur allgemeinen und deskriptiven Mengenlehre, zur Maßtheorie und zur Analysis geleistet. Auch seine Beiträge zur Theorie der Lie-Algebren, zur Wahrscheinlichkeitstheorie und zur Versicherungsmathematik sind für die Folgeentwicklung bedeutsam gewesen. Hausdorff war aber auch – für einen Mathematiker ziemlich einmalig – ein bemerkenswerter Literat und philosophischer Schriftsteller. Von 1897 bis 1913 publizierte er unter dem Pseudonym Paul Mongré einen Aphorismenband, ein erkenntniskritisches Buch, einen Gedichtband, ein in über 30 Städten mehr als 300 mal aufgeführtes Theaterstück sowie 17 Essays in damals führenden literarischen Zeitschriften. In den zwanziger und frühen dreißiger Jahren war er – weltweit anerkannt und geachtet – der führende Repräsentant des Bonner Mathematischen Instituts. Als Jude unter der nationalsozialistischen Diktatur verfolgt und gedemütigt, nahm er sich gemeinsam mit seiner Frau am 26. Januar 1942 das Leben, als die Deportation in ein Konzentrationslager unmittelbar bevorstand.

    Egbert Brieskorns Engagement begann 1979, als aus der Studentenschaft des Mathematischen Instituts der Vorschlag kam, Felix Hausdorff durch eine Gedenktafel im Institut zu ehren. Brieskorn, der sich schon viele Jahre in der Friedensbewegung engagiert und sich eingehend mit Fragen der Verantwortung des Wissenschaftlers in der Gesellschaft auseinandergesetzt hatte, unterstützte das Vorhaben von Beginn an, leistete selbst finanzielle Hilfe und organisierte auch im Lehrkörper eine Spendensammlung. So konnte 1980, zu Hausdorffs Todestag am 26. Januar, die marmorne Tafel eingeweiht werden, die sich über 30 Jahre im alten Institut in der Wegelerstraße befand und die kürzlich im neuen Institutsgebäude angebracht wurde. Aus Anlass der Einweihung erschien ein Artikel über Hausdorffs tragisches Schicksal aus der Feder des Historikers Herbert Mehrtens. Egbert Brieskorn verfasste dazu ein Vorwort, dessen Schlusssätze hier zitiert seien; er schreibt dort: "`Keine Form von Unmenschlichkeit und Unterdrückung kann uns gleichgültig lassen, auch wenn die Opfer uns ferne und unbekannt sind. Der Gedanke an einen Menschen wie Hausdorff, den wir alle wegen seiner großen wissenschaftlichen Leistung kennen, kann dazu beitragen, unser Gewissen und unsere Verantwortlichkeit zu schärfen. Ohne ein wachsendes Gefühl der Verantwortung werden die Wissenschaftler nicht in der Lage sein, ihren Beitrag zu einer menschlicheren Zukunft zu leisten. Für diese Aufgabe brauchen wir auch die 'Bewältigung' der Vergangenheit."'

      Im November 1980 konnte Egbert Brieskorn es erreichen, dass der Nachlass Hausdorffs, der einen Umfang von etwa 26.000 Blatt hat und sich seit 1964 im Privatbesitz von Prof. Günther Bergmann in Münster befand, an die Universitätsbibliothek Bonn verkauft wurde. Der Erlös kam Hausdorffs Tochter Lenore König, die in ärmlichen Verhältnissen in Bonn in einem Altenheim lebte, zugute. Als der Vertrag unterschrieben war, schrieb Egbert Brieskorn am 15. November an Günther Bergmann: "`Die Regelungen, die jetzt getroffen worden sind, scheinen mir sehr gut. Ich bin zwar selber kein Historiker, habe aber doch ein gewisses Interesse für die Geschichte der Mathematik und vergebe gelegentlich auch Diplomarbeiten und Dissertationen mit historischen Aspekten. Ich hoffe, dass sich eines Tages auch ein Mathematikhistoriker findet, der über Hausdorff arbeiten will. Dann wird der Nachlass in unserer Universitätsbibliothek sehr wichtig werden."' Er hat bei diesen Zeilen gewiss nicht daran gedacht, dass er es selbst sein würde, der fast 10 Jahre später das Projekt in Angriff nehmen würde, eine Biographie Hausdorffs zu schreiben.  

     Am 26. Januar 1992 jährte sich Hausdorffs Todestag zum 50. Male. Aus diesem Anlass richtete Egbert Brieskorn ein Gedenkkolloquium des Mathematischen Instituts der Universität Bonn aus, als dessen Ergebnis der von ihm herausgegebene Band "`Felix Hausdorff zum Gedächtnis – Aspekte seines Werkes"' im Vieweg-Verlag erschien. Ferner gestaltete er eine Ausstellung über Leben und Werk Felix Hausdorffs, die über das Mathematische Institut hinaus an der Universität und in der Bonner Öffentlichkeit lebhaftes Interesse fand. Zur Ausstellung publizierte er einen umfangreichen Katalog mit einer ersten biographischen Skizze Hausdorffs. Unter dem Titel "`Auf dünner Schneide tanzt mein Glück"', der ersten Zeile von Hausdorffs Gedicht "`Wiederkunft“ aus dem Gedichtband "`Ekstasen"', gab es am 1. Februar 1992 in der Sendereihe "`Mosaik"' des WDR eine einstündige Sendung über Hausdorff und die Bonner Ausstellung, konzipiert von Egbert Brieskorn und dem Redakteur für Kultur Friedrich Riehl. Sie hat die Ausstellung und ihr Anliegen weit über Bonn hinaus bekannt gemacht. Auch das Presseecho auf die Ausstellung war beachtlich.

   Die Vorbereitungen für alle diese Aktivitäten begannen bereits 1989 mit zahlreichen Gesprächen, die Egbert Brieskorn mit Hausdorffs Tochter führte, mit Kontakten zu Hausdorffs Nichte Else Pappenheim und zu weiteren Zeitzeugen sowie mit dem Sammeln von Materialien zur Biographie, insbesondere mit Recherchen in Archiven. Bei diesen Recherchen, die er über viele Jahre fortgesetzt hat, hat er eine Gründlichkeit und einen detektivischen Spürsinn bewiesen, der jedem professionellen Historiker alle Ehre machen würde. In seinem Nachlass gibt es Dutzende von dickleibigen Ordnern, die alle diese Bemühungen dokumentieren und im Erfolgsfall, der ja durchaus nicht immer eintritt, ihre Ergebnisse festhalten.

  Ebenfalls im Vorfeld des Gedenkkolloquiums haben die Bonner Mathematiker Egbert Brieskorn, Friedrich Hirzebruch und Stefan Hildebrandt unter Hinzuziehung auswärtiger Experten über Möglichkeiten und notwendige Schritte beraten, eine Edition der Werke Hausdorffs in die Wege zu leiten. Friedrich Hirzebruch schlug vor, bei der Nordrhein-Westfälischen Akademie der Wissenschaften, die bereits eine Reihe von Editionsprojekten betreute, eine Hausdorff-Kommission zu bilden, die ein solches Projekt planen und dann auch betreuen sollte. Die Akademie stimmte zu und die Kommission nahm unter der Leitung von Akademiemitglied Reinhold Remmert noch 1991 ihre Arbeit auf.  Als erster Schritt wurde eine sorgfältige inhaltliche Erschließung und Katalogisierung des Hausdorffschen Nachlasses ins Auge gefaßt. Diese erfolgte von Oktober 1993 bis Ende 1995. Das Ergebnis war ein Findbuch von 550 Seiten mit kurzen Beschreibungen des Inhalts jedes einzelnen Faszikels. Nachdem diese Voraussetzung für ein erfolgreiches Editionsprojekt unter Einbeziehung des Nachlasses geschaffen war, konnte man einem solchen Vorhaben näher treten.

   Um eine sorgfältig kommentierte Edition mit Einschluss ausgewählter Teile aus dem Nachlass ins Werk zu setzen, muss vieles bedacht und getan werden: Man muss Editionsprinzipien aufstellen, eine Bandstruktur erarbeiten, geeignete Mitarbeiter suchen und für die Sache gewinnen, man muss Anträge stellen, um die Sache zu finanzieren usw. usf. Egbert Brieskorn war der spiritus rector bei all diesen Debatten und Aktivitäten. Besonders schwierig gestaltete sich die Gewinnung geeigneter Philosophen und Literaturwissenschaftler als Editoren für die Bände, welche dieser Seite von Hausdorffs Schaffen gewidmet sind. Er nahm zu diesem Zweck an einer philosophischen Tagung zum jüdischen Nietzscheanismus in Greifswald teil, um mit entsprechenden Experten ins Gespräch zu kommen, und organisierte in Bonn mit Unterstützung des MPI für Mathematik eine interdisziplinäre Tagung mit Geisteswissenschaftlern über Hausdorffs philosophisches und literarisches Werk. 
 
    Der Antrag auf Finanzierung, den er 1996 gemeinsam mit Friedrich Hirzebruch und Erhard Scholz bei der DFG stellte, wurde schließlich genehmigt und im November 1996 nahm die Arbeitsstelle "`Hausdorff-Edition"' am Mathematischen Institut unter seiner Leitung die Arbeit auf. 2002 übernahm die Nordrhein-Westfälische Akademie der Wissenschaften und der Künste die Hausdorff-Edition als eines ihrer Akademieprojekte. Einige der ursprünglich gewonnenen Editoren haben aus verschiedenen Gründen dann doch nicht mitarbeiten können, so dass im Laufe der Zeit auch immer wieder einmal ein neuer Mitarbeiter gefunden werden musste. Schließlich haben an der Edition vierzehn Mathematiker, vier Mathematikhistoriker, zwei Literaturwissenschaftler, ein Philosoph und ein Astronom aus vier Ländern mitgewirkt.

      Es ging Egbert Brieskorn bei der Edition besonders darum, den vielen sichtbaren aber auch den auf den ersten Blick unsichtbaren Fäden nachzuspüren, die von Hausdorff als Mathematiker zum Literaten und Philosophen Mongré hinüber- und herüberführen. Er betonte also immer wieder den interdisziplinären Charakter des Projektes, hielt auf drei großen Editorenkonferenzen wegweisende Vorträge gerade in dieser Hinsicht und organisierte eine Reihe von Diskussionsrunden mit beteiligten Mathematikern, Philosophen und Literaturwissenschaftlern. Um einen Eindruck von seinen Intentionen zu gewinnen, zitieren wir hier den Beginn seines programmatischen Vortrags auf einer dieser Konferenzen im Februar 2003 auf Schloß Rauischholzhausen: "`Es ist nicht immer ausgemacht, dass der Gedanke in den Zeilen oder auch nur zwischen den Zeilen stehen müsse; vielleicht steht er ganz woanders, weit, weit entfernt! Vielleicht zieht der Autor seine Glocke, damit irgendwo eine Saite in gleiche Schwingungen gerathe und in ihrer Klangfarbe antworte – und nicht der Glockenton, sondern der Saitenklang spricht den eigentlichen Gedanken aus."'  Diesen Aphorismus von Paul Mongré sollten wir bei unserem Rundgespräch über diesen Autor, über den Schriftsteller  Mongré und über den Mathematiker Felix Hausdorff, immer im Sinn behalten. Jeder von uns wird aus der Vielfalt der Motive und dem Reichtum der Klangfarben etwas anderes heraushören. Was man hört, hängt ja sehr von den Hörerfahrungen ab, die man im Laufe seines Lebens gemacht hat. Wenn Literaturwissenschaftler, Philosophen, Historiker und Mathematiker auch aufeinander hören, können wir hoffen, hier und da den 'eigentlichen Gedanken' herauszuhören."'

     Als Leiter der Arbeitsstelle hat Egbert Brieskorn dem Koordinator des Editionsprojektes und den Editoren und Mitarbeitern der einzelnen Bände alle Freiheiten gelassen. Wenn Schwierigkeiten aufgetreten sind, hat er mit Rat und Tat geholfen. Einige eingereichte Beiträge hat er aber auch zurückgewiesen, weil sie seinen hohen Anforderungen nicht genügten. Es ist jedoch in diesen Fällen stets gelungen, entscheidende Verbesserungen zu erzielen. Mittlerweile sind acht der zehn geplanten Bände bei Springer erschienen: Band IA "`Allgemeine Mengenlehre"' (2013), Band II "`Grundzüge der Mengenlehre"' (2002), Band III "`Deskriptive Mengenlehre und Topologie"' (2008), Band IV "`Analysis, Algebra und Zahlentheorie"' (2001), Band V "`Astronomie, Optik und Wahrscheinlichkeitstheorie"' (2006), Band VII "`Philosophisches Werk"' (2004), Band VIII "`Literarisches Werk"' (2010), Band IX "`Korrespondenz"' (2012). Egbert Brieskorns Handschrift ist in allen Bänden spürbar, auch dann, wenn er nicht als Band-Mitarbeiter in Erscheinung getreten ist. 
   
      Den schwierigsten Teil des ganzen Projektes aber hat er sich selbst aufgeladen: Band IB, die Biographie. Hier waren neben der mathematischen Arbeit Hausdorffs auch weite Felder seiner Interessen und lebensweltlichen Bezüge aus sehr verschiedenen Bereichen zu berücksichtigen: Philosophie, vor allem Kant, Schopenhauer, Nietzsche und Hausdorffs Beziehungen zum Nietzsche-Archiv, Erkenntniskritik, vor allem Hausdorffs Sprachkritik und seine Überlegungen zu Raum und Zeit, Hausdorffs literarisches Werk und seine Beziehungen zu Literaten wie Dehmel, Hartleben und Wedekind, Musik, insbesondere Hausdorffs Verhältnis zu Wagner und seine Beziehung zu Reger, bildende Kunst, insbesondere Hausdorffs Freundschaft mit Max Klinger. Hinzu kommen die Familiengeschichte im Rahmen der jüdischen Geschichte und die Geschichte des Antisemitismus bis zu Hausdorffs tragischem Ende unter der Nazidiktatur. 

    2007 hatten Egbert Brieskorn, Erhard Scholz und der Koautor dieses Nachrufs die Gelegenheit, im Séminaire d'Histoire des Mathématiques de l'Institut Henri Poincaré in Paris das Projekt der Hausdorff-Edition vorzustellen. Dort hielt Brieskorn den einleitenden Vortrag, in dem er rückschauend auf sein eigenes Wirken folgendes sagte (er sprach natürlich französisch; wir zitieren aus seinem deutschsprachigem Entwurf): "`Was nun meinen eigenen Anteil an dem Projekt angeht, so will ich gleich gestehen, dass ich kein Historiker bin und dass mich nicht vorwiegend ein historisches Interesse leitet. Vielmehr entsprang mein persönliches Interesse ursprünglich zwei Motiven. Das eine Motiv war die Scham über die entsetzliche Schuld, die Deutschland mit der Verfolgung und Ermordung der Juden in Europa auf sich geladen hat. Das andere Motiv war ganz persönlich: In den achtziger Jahren lernte ich Felix Hausdorffs Tochter Lenore König kennen, die damals in einem Altersheim in Bonn lebte. Sie öffnete mir einen ersten Zugang zum Leben und zur Persönlichkeit ihres Vaters. Dadurch empfand ich später die persönliche Verpflichtung, das Leben dieses ungewöhnlichen Menschen besser zu verstehen. Als die Universität Bonn eine Gedenkveranstaltung zu Hausdorffs fünfzigstem Todestag vorbereitete, und später, als der Plan einer Edition seiner Werke Gestalt annahm, sah ich die Biographie Felix Hausdorffs als meine persönliche Aufgabe an. Ich habe diese Aufgabe anfangs unterschätzt, sowohl hinsichtlich der Schwierigkeit wie auch im Hinblick darauf, was diese Arbeit für mich persönlich bedeutete. Diese Arbeit hat auch mein eigenes Leben und Denken sehr verändert: Ich habe – so hoffe ich – manches von Hausdorff gelernt."'

     An der Niederschrift der Biographie hat Egbert Brieskorn die letzten Jahre intensiv gearbeitet, auch und gerade während seiner schweren Krankheit. Drei Wochen vor seinem Tod hat er der Bonner Arbeitsstelle das letzte Unterkapitel geschickt, das er noch fertigstellen konnte, übrigens ein besonders schwieriges, über Hausdorffs Beziehungen zu dem Philosophen, Mystiker und Anarchisten Gustav Landauer. Bis jetzt liegen aus seiner Feder 530 Seiten der Biographie druckreif geschrieben vor. Als er spürte, dass er mehr nicht schaffen würde, hat er ein Gespräch vorgeschlagen, um zu erzählen, wie er sich den weiteren Verlauf der Biographie vorgestellt hat. Als Termin hatte er den 12. Juli 2013 vorgeschlagen. Am Abend des 11. Juli ist er verstorben. Es ist für die Bonner Arbeitsstelle eine Verpflichtung, den Band zu Ende zu bringen, so gut wir es vermögen. Dabei werden die über 100 Ordner mit Material zur Biographie, die er in mehr als 20 Jahren in oft mühevoller Kleinarbeit gesammelt hat und die uns Frau Brieskorn übergeben hat, eine unschätzbare Hilfe sein.
\bigskip

\textbf{\large Brieskorns mathematisches Werk}

Das mathematische Werk von Brieskorn ist weitgehend bestimmt durch die Entwicklung der \textit{Singularitätentheorie komplexer Hyperflächen}. Insbesondere seine frühen Arbeiten hatten einen großen Einfluss auf die Entwicklung der Singularitätentheorie und es ist keine Übertreibung, Brieskorn zusammen mit Vladimir Igorevich Arnold, John W. Milnor und René Thom als einen der Väter der Singularitätentheorie zu bezeichnen. In seinem Übersichtsvortrag \cite{EB1976} schreibt Brieskorn:

\textit{Singularitäten gibt es in allen Gebieten der Mathematik und in vielen Anwendungen, und das Gegensatzpaar "`regulär - singulär"' ist eines der häufigsten in einer ganzen Reihe von solchen Gegensatzpaaren in der Mathematik. Was mit Singularitäten eigentlich gemeint ist, zeigt die Analyse der vielen verschiedenen Definitionen singulärer Objekte. Eine solche Analyse führt auf einige wenige Grundbedeutungen: Eine Singularität innerhalb einer Gesamtheit ist eine Stelle der Einzigartigkeit, der Besonderheit, der Entartung, der Unbestimmtheit oder der Unendlichkeit. Alle diese Bedeutungen hängen eng miteinander zusammen.}

Den Begriff "`Singularitätentheorie"' verwende ich hier im Sinne der Untersuchung von Systemen endlich vieler differenzierbarer, analytischer oder algebraischer Funktionen in der Umgebung eines Punktes, in dem die Jacobimatrix der Funktionen nicht maximalen Rang hat. Betrachtet man die Nullstellenmenge der Funktionen, so bedeutet dies nach dem Satz über implizite Funktionen, dass diese in der Nähe eines singulären Punktes keine differenzierbare, analytische oder algebraische Mannigfaltigkeit ist. Hierbei sind Singularitäten von Vektorfeldern oder Differentialformen mit eingeschlossen.

Auch wenn ich den Begriff Singularitätentheorie gebrauche, der meines Wissens von V.I. Arnold eingeführt wurde, so handelt es sich nicht wirklich um eine geschlossene Theorie mit mehr oder weniger einheitlichen Methoden. Im Gegenteil, ein Charakteristikum dieses Gebietes ist die Vielfalt der verwendeten Methoden und der Beziehungen zu anderen mathematischen Disziplinen. Dazu gehören die algebraische Geometrie, komplexe Analysis, kommutative Algebra, Kombinatorik, Darstellungstheorie, Liegruppen, Topologie, Differentialtopologie, dynamische Systeme, symplektische Geometrie und andere. Gerade die Vielschichtigkeit und die vielfältigen Wechselwirkungen der Singularitätentheorie mit anderen mathematischen und nicht-mathematischen Gebieten, wie z.B. der geometrischen Optik, haben Brieskorn besonders fasziniert und er hat, wie wir sehen werden, zu der Erforschung einiger dieser Wechselwirkungen wesentlich beigetragen. Mit dem Begriff Singularitäten{\em theorie} hat er sich allerdings nie anfreunden können.

Fast alle mathematischen Arbeiten von Brieskorn beschäftigen sich entweder direkt mit Singularitäten komplexer Hyperflächen, oder sie sind durch das Studium dieser Singularitäten motiviert. Dabei zeigt sich in seinem Werk neben Originalität und Tiefe eine große Breite der Fragestellungen und der verwendeten Methoden, wie sie für das gesamte Gebiet typisch ist.

Bei der folgenden Übersicht über Brieskorns Arbeiten versuche ich auch wichtige Ergebnisse, auf denen Brieskorns Resultate beruhen, sowie Entwicklungen, die aus seinen Arbeiten entstanden, mit aufzuzeigen. Eine kurze Darstellung des wissenschaftlichen Werkes von Brieskorn findet sich in \cite{GG1997}. 
\bigskip
\medskip

\textbf{Dissertation}

Singularitäten spielen in den ersten beiden Arbeiten von Brieskorn allerdings noch keine Rolle. Es handelt sich dabei um Teile seiner Dissertation \cite{EB1962}, die er 1962 als Student von Friedrich Hirzebruch verfasste und die er in \cite{EB1964} und \cite{EB1965} veröffentlichte.

Die Fragestellung dort ist, inwieweit die einer Kähler--Mannigfaltigkeit zugrunde liegende differenzierbare Strukur im Fall der komplexen Quadriken $Q_n$ bzw. der holomorphen $\P^n$--Bündel über $\P^1$ bereits die biholomorphe Stuktur bestimmt. Dieses Problem war von F. Hirzebruch und K. Kodaira im Jahre 1957 für den komplex projektiven Raum $\P^n$  untersucht und gelöst worden.

Das folgende Hauptergebnis der ersten Arbeit ist ein genaues Analogon des erwähnten Satzes von Hirzebruch und Kodaira. Das Ergebnis hatte Brieskorn bereits 1961 in den Notices der AMS angekündigt:

\textit{Sei $X$ eine $n$--dimensionale Kählersche Mannigfaltigkeit, die diffeomorph zur $n$--dimensionalen komplex projektiven Quadrik $Q_n$ ist. Dann gilt:}
\textit{
\begin{enumerate}
\item [(i)] Ist $n$ ungerade, so ist $X$ biholomorph zu $Q_n$
\item [(ii)] Ist $n$ gerade und $n\neq 2$, dann gilt für die 1. Chernsche Klasse $c_1$ von $X$: $c_1=\pm ng$, wobei $g$ das positive Erzeugende von $H^2(X, \Z)\cong \Z$ bezeichnet; falls $c_1=+ng$ gilt, dann ist $X$ biholomorph zu $Q_n$.
\end{enumerate}}

Brieskorn fragt, ob die Voraussetzung, dass $X$ Kählersch ist, weggelassen werden kann und ob es überhaupt Kähler--Mannigfaltigkeiten mit $c_1=-ng$ gibt, die diffeomorph zu $Q_n$ sind. Beide Probleme scheinen bis heute offen zu sein.

Als Folgerung beweist Brieskorn Aussagen über das Deformationsverhalten von $Q_n$. Er betrachtet eine Familie komplexer Mannigfaltigkeiten $V_t, t\in \C,\ |t|$ hinreichend klein, und beweist: 
\textit{
\begin{enumerate}
\item [(i)] Ist $V_0\cong Q_n$, dann ist auch $V_t\cong Q_n$ für $t\neq 0$,
\item [(ii)] Ist $V_0$ Kählersch und $V_t\cong Q_n$ für $t\neq 0$ und $n\geq 3$, dann ist $V_0\cong Q_n$.
\end{enumerate}}

Auch hier stellt er die Frage, ob in (ii) die Voraussetzung, dass $V_0$ Kählersch ist, weggelassen werden kann.

Dass dies tatsächlich der Fall ist, wurde 1995 von J. M. Hwang bewiesen, nachdem dieselbe Frage der "`Nichtdeformierbarkeit"' von $\P^n$ statt $Q_n$ bereits vorher von Y.-T. Siu positiv beantwortet worden war.

Für $n=2$ gilt die erste Aussage nicht, denn auf der differenzierbaren Mannigfaltigkeit, die $Q_2$ zugrunde liegt, gibt es nach Hirzebruch unendlich viele verschiedene komplexe Strukturen, die sogenannten Hirzebruchschen $\Sigma$--Flächen $\Sigma_{2m}$.

Die $\Sigma$--Flächen sind Totalräume von holomorphen Faserbündeln über $\P^1$ mit Faser $\P^1$. Im zweiten Teil seiner Dissertation untersucht Brieskorn holomorphe Faserbündel über $\P^1$ mit Faser $\P^n$, die er in Anlehnung an die Hirzebruchschen $\Sigma$--Flächen \textit{$\Sigma$--Mannigfaltigkeiten} nennt. Unter Ausnutzung des Spaltungssatzes von Grothendieck für Vektorbündel über $\P^1$ klassifiziert Brieskorn dort die $\Sigma$--Mannigfaltigkeiten bis auf biholomorphe und birationale Äquivalenz und bis auf Diffeomorphie. Daraus ergibt sich, dass es wie im Fall von $\Sigma$--Flächen auf jeder differenzierbaren $\Sigma$--Mannigfaltigkeit abzählbar unendlich viele verschiedene komplexe Strukturen gibt. Außerdem beweist er, dass $\Sigma$--Mannigfaltigkeiten wieder in $\Sigma$--Mannigfaltigkeiten deformieren und dass in einer Kählerschen Familie von $\Sigma$--Flächen diese in eine $\Sigma$--Fläche spezialisieren (ähnlich wie in (ii) oben).

Die Fragestellungen, wie die Beweismethoden der Dissertation, stammen aus der algebraischen und analytischen Geometrie sowie der algebraischen Topologie. Diese Methoden, einschließlich der aus Frankreich stammenden Garbentheorie, waren damals ganz neu und begannen sich in Deutschland erst langsam durchzusetzen, vor allem in der Generation junger Mathematiker. In exemplarischer Weise wurden sie von Brieskorns Lehrer Friedrich Hirzebruch verkörpert. Neben Hirzebruch hatten auch Hans Grauert und Reinhold Remmert mit der Entwicklung der Theorie allgemeiner komplexer Räume, deren Strukturgarbe nilpotente Elemente enthalten konnte, großen Einfluss auf die Entwicklung der modernen analytischen und algebraischen Geometrie in Deutschland. Brieskorn, der sich in seiner Dissertation neben Hirzebruch auch bei Reinhold Remmert und Antonius van de Ven bedankt, die er beide während seiner Tätigkeit als wissenschaftliche Hilfskraft und Mitarbeiter 1962 in Erlangen traf, hatte also ein äußerst inspirierendes Umfeld moderner Mathematik. Er wurde von der Aufbruchstimmung, die damals in Deutschland herrschte, und besonders durch die Förderung seines Lehrers Hirzebruch entscheidend beeinflusst.
\medskip

\textbf{Deformationstheorie}

Friedrich Hirzebruchs Buch {\em "`Neue topologische Methoden in der algebraischen Geometrie"'}, das bereits 1956 in der Springer Reihe "`Ergebnisse der Mathematik und ihrer Grenzgebiete"' veröffentlicht worden war, war gerade in der zweiten, erweiterten Auflage erschienen und der dort bewiesene Satz von Hirzebruch--Riemann--Roch bildete eine der Grundlagen von Brieskorns Dissertation. Eine weitere Grundlage bildete die von K. Kodaira und D.C. Spencer entwickelte Deformationstheorie analytischer Strukturen.

Der Satz von Hirzebruch--Riemann--Roch war eine großartige Verallgemeinerung des klassischen Satzes von Riemann--Roch auf komplexe Vektorbündel auf beliebigen komplex projektiven Mannigfaltigkeiten, statt Divisoren auf Riemannschen Flächen, mit den Methoden der sich damals gerade durchsetzenden Garbentheorie. Wie Brieskorn in einem Lebenslauf schreibt, war dieser Satz der Grund, warum er von München nach Bonn zu Hirzebruch wechselte. 1963 verallgemeinerten M. Atiyah und I. Singer den Satz von Hirzebruch--Riemann--Roch  zu einem Indexsatz über elliptische Differentialoperatoren auf einer komplexen Mannigfaltigkeit, der bedeutende Sätze der Differentialgeometrie umfasst und wichtige Anwendungen in der theoretischen Physik hat, wofür sie 2004 den Abelpreis erhielten. Dieser Satz wurde 1967 von Grothendieck algebraisiert und in funktorieller Weise auf eigentliche Morphismen quasiprojektiver Schemata erweitert.  Modifikationen bzw. Weiterentwicklungen halten bis heute an, z.B. mit dem "`Quantum Riemann--Roch"' in der Gromov--Witten--Theorie.  

Die von Kodaira und Spencer in \cite{KS1958} entwickelte Deformationstheorie komplexer Mannigfaltigkeiten, insbesondere die infinitesimale Theorie der Linearisierung von Deformationen sowie der Satz von Kuranishi über die Existenz einer semi-universellen Deformation, wurden später weiter entwickelt und gehören heute zu den wichtigsten Methoden in der komplexen Analysis sowie der algebraischen und arithmetischen Geometrie. Zu erwähnen sei hier nur die von Hans Grauert und Adrien Douady bewiesene Existenz einer semi-universellen Deformation für kompakte komplexe Räume.

Noch wichtiger für die Singularitätentheorie und für Brieskorns spätere Arbeiten war der von Grauert bewiesene Satz über die Existenz einer semi-universellen Deformation isolierter Singularitäten komplexer Räume. An der Entstehung des Beweises durch Grauert war Brieskorn selbst beteiligt. Und das kam wie folgt.

Brieskorn hatte im Juli 1969 seine Professur in Göttingen angetreten, war aber von September 1969 bis Februar 1970 zu einem Forschungsaufenthalt an das IHES in Bures-sur-Yvette beurlaubt worden. Von dort brachte er ein interessantes Problem für das gemeinsam von Brieskorn und Grauert geleitete Oberseminar im Sommer- und Wintersemester 1970/71 nach Göttingen mit: den Nachweis der Existenz einer semi--universellen Deformation isolierter Singularitäten komplexer Räume. Als Grundlage sollten die Arbeiten von M. Schlessinger "`Functors of Artin Rings"' und von G. N. Tyurina "`Locally semi--universal flat deformations of isolated singularities of complex spaces"' dienen. Tyurinas Arbeit war 1969  auf russisch, aber erst 1971 in englischer Übersetzung erschienen. Der Verfasser dieser Zeilen, der selbst an dem Oberseminar teilnahm, glaubt sich zu erinnern, dass Brieskorn beide in Göttingen noch nicht bekannte Arbeiten mitgebracht und insbesondere eine englische Übersetzung von Tyurinas Arbeit beschafft hatte. 

Schlessinger hatte in seiner Arbeit Bedingungen für die Existenz einer formalen semi-universellen Deformation angegeben, die "`Schlessinger--Bedingungen"', während Tyurina die Existenz für normale isolierte Singularitäten bewiesen hatte, unter der Zusatzbedingung, dass die zweite Ext--Gruppe der holomorphen $1$--Formen auf der Singularität verschwindet. Tyurina hatte die Arbeit von Schlessinger offenbar nicht gekannt und Brieskorns Idee war es, dass die Verbindung der Ansätze von Schlessinger und Tyurina den Beweis für beliebige isolierte Singularitäten ohne Tyurinas Zusatzbedingungen liefern sollte. Er hatte sogar schon genauere Vorstellungen, wie beide Arbeiten zu einem Beweis zusammengeführt werden sollten und verteilte die Vorträge entsprechend an die Teilnehmer des Oberseminars. Der letzte Vortrag ging an Grauert mit Brieskorns Hinweis "`und Sie beweisen dann den allgemeinen Satz ohne Voraussetzung"'. Grauert antwortete nur "`dafür brauche ich aber die Weihnachtsferien"', was allgemeine Heiterkeit auslöste. Die Mitarbeit in dem Oberseminar war ziemlich anspruchsvoll, besonders für jemanden, der gerade erst mit der Arbeit an seiner Diplomarbeit begonnen hatte, aber gleichzeitig enorm anregend und bereichernd. Allen Teilnehmern war bewusst, dass sie an der Entstehung eines bedeutenden Resultates beteiligt waren und warteten gespannt auf Grauerts Vortrag, der Anfang 1971 stattfinden sollte. Grauert begann seinen Vortrag dann aber mit der Bemerkung, dass er den Beweis leider nicht vortragen könne und er auch nicht sicher sei, ob man vielleicht doch eine Voraussetzung wie die von Tyurina brauche. Allerdings könne er eine interessante Verallgemeinerung des Weierstraßschen Divisionssatzes vortragen, die Division mit Rest nach einem Ideal (ein Ergebnis, das unabhängig von Hironaka im Zusammenhang mit der Auflösung von Singularitäten gefunden war). Dieser allgemeine Divisionssatz bildete dann das entscheidende Hilfsmittel für den Beweis der Existenz einer semi-universellen Deformation isolierter Singularitäten ohne jede Voraussetzung, den Grauert im Sommer 1971 vollendete (veröffentlicht in \cite{HG1972}). 

Die Episode zeigt das untrügliche Gespür Brieskorns für interessante und wichtige mathematische Probleme, was man in seinem ganzen Werk und auch bei der Auswahl der Themen für Diplom-- und Doktorarbeiten beobachten kann.
\medskip

\textbf{Quotientensingularitäten und simultane Auflösung}

Die Jahre nach der Promotion gehörten zu den fruchtbarsten in Brieskorns wissenschaftlichem Leben. Irgendwann im Jahr 1963 hatte Hirzebruch Brieskorn vorgeschlagen, die Arbeit "`On analytic surfaces with double points"' von Michael Atiyah \cite{MA1958} zu studieren und zu verallgemeinern. In dieser Arbeit hatte Atiyah unter anderem gezeigt, dass eine Familie $f\!\!:\!\!X\to S$ kompakter komplexer Flächen über einer glatten $1$--dimensio\-nalen Mannigfaltigkeit $S$, deren allgemeine Faser glatt ist und deren endlich viele spezielle Fasern nur Singularitäten vom Typ $A_1$ haben, eine simultane Auflösung besitzt. Hierbei ist eine simultane Auflösung einer allgemeinen holomorphen Abbildung $f\!\!:\!\!X\to S$ ein kommutatives Diagramm holomorpher Abbildungen
\[
\xymatrix{
Y\ar[r]^\psi\ar[d]_g & X\ar[d]^f\\
T\ar[r]^\varphi & S\ ,}
\]
wobei $\varphi$ eine verzweigte Überlagerung und $g$ eine nicht--singuläre,  eigentliche surjektive Abbildung ist, die für alle Fasern $X_s\!\!=\!\!f^{-1}(s)$ von $f$ eine Auflösung $\psi|Y_t\!\!: Y_t\to X_s\ ,\ \varphi(t)=s$, induziert.

Wegen der lokalen Monodromie um die singulären Fasern von $f$ ist der Basiswechsel $\varphi$ nötig, d.h. eine simultane Auflösung von $f$ über $S$ selbst kann es i.A. nicht geben. Die lokale Monodromie um eine singuläre Faser $X_{s_0}$ sieht man so: schränkt man $f$ auf einen kleinen geschlossenen Weg $\gamma$ mit Anfangs- und Endpunkt $s$ um $s_0\in S$ ein, so dass $\gamma$ durch keinen singulären Wert von $f$ geht und $s_0$ einfach umrundet, so erhält man ein lokal triviales Faserbündel über $\gamma$, das i.A. nicht trivial ist. Durch Wegliftung (z.B. mittels eines Ehresmannschen Zusammenhangs) erhält man einen nicht--trivialen Diffeomorphismus der Faser $X_s$, die \textit{ geometrische Monodromie}, die einen nicht--trivialen Isomorphismus der Homologie von $X_s$  induziert. Der Basiswechsel $\varphi$ eliminiert die lokale Monodromie, denn da alle Fasern von $g$ nicht--singulär sind, ist die Monodromie von $g$ über $T$ trivial.

Hirzebruch hatte Brieskorn vorgeschlagen, die Arbeit von Atiyah auf Familien von Flächen mit Singularitäten vom Type $A_k, D_k, E_6, E_7, E_8$ zu verallgemeinern. Es stellte sich heraus, dass dies eine wunderbare Idee war, und es begann eine hochinteressante Geschichte mit vielen Akteuren und großartigen Entdeckungen. Am Ende stand nicht nur die simultane Auflösung der Familien von Flächen mit $ADE$--Singularitäten sondern auch die Entdeckung von exotischen Sphären als Umgebungsränder von Singularitäten. Brieskorn schreibt in \cite{EB2010}  \textit{"`I shall be grateful for it to my teacher till the day that I die"'}.

Da das Problem lokal ist, genügt es, eine Abbildung $f$ von $X=\C^3$ nach $S=\C$ (bzw. von kleinen Umgebungen der Nullpunkte) der Form $s=f(x,y,z)$ zu untersuchen. Hierbei ist $f(x,y,z)=0$ die Gleichung einer \textit{$ADE$--Singularität}, $f$ also ein Polynom der folgenden Liste:
\[
\begin{array}{lll}
A_k: & x^{k+1}+y^2+z^2, & k\geq 1\\
D_k: & x^{k-1}+xy^2+z^2, & k\geq 4\\
E_6: & x^4+y^3+z^2\\
E_7: & x^3y+y^3+z^2,\\
E_8: & x^5+y^3+z^2\ .
\end{array}
\]

Diese Polynome tauchen schon im 19. Jahrhundert in den Arbeiten von Hermann Amandus Schwarz und Felix Klein auf (vgl. \cite{FK1884}). Seit Kleins Zeiten erscheinen sie in immer neuen, verschiedenen Zusammenhängen (vgl. \cite{GG1992} für einen Überblick) und haben dadurch die Mathematiker bis heute fasziniert. Je nachdem, in welchem Zusammenhang sie betrachtet werden, heißen sie auch \textit{einfache Flächensingularitäten}, \textit{rationale Doppelpunkte}, \textit{Du-Val-Singularitäten} oder \textit{Kleinsche Singularitäten}.

Der Zusammenhang, in dem die $ADE$--Singularitäten bei Klein auftauchen, ist für uns besonders interessant. Klein klassifizierte die endlichen Untergruppen $G$ von $\SL(2, \C)$ bis auf Konjugation und erhielt die folgenden Gruppen:
\[
\begin{array}{ll}
C_{k+1}: & \text{die zyklische Gruppe der Ordnung } k+1\\
D_{k-2}: & \text{die binäre Diedergruppe der Ordnung } $4(k-2)$\\
T: & \text{die binäre Tetraedergruppe der Ordnung } $24$\\
O: & \text{die binäre Oktaedergruppe der Ordnung } $48$\\
I: & \text{die binäre Ikosaedergruppe der Ordnung } $120$ 
\end{array}
\]

Diese Gruppen sind (komplexe) \textit{Spiegelungsgruppen}, also von komplexen Spiegelungen (Automorphismen endlicher Ordnung, die eine Hyperebene festhalten) erzeugte Gruppen, und Klein bewies, dass der Ring $\C[z_1, z_2]^G$ der unter $G$ invarianten Polynome in $\C[z_1, z_2]$ durch drei invariante Polynome $X,Y,Z$ erzeugt wird, die genau einer Relation $f(X,Y,Z)=0$ genügen. Klein bestimmte die Relationen und fand, dass diese für die Gruppen $C_{k+1}, D_{k-2}, T, O, I$ durch die Polynome $A_k, D_k, E_6, E_7, E_8$ gegeben sind.

Das bedeutet, dass die ADE--Singularitäten \textit{$2$--dimensionale Quotientensingularitäten} sind: Sei allgemeiner $G\subset \GL(2, \C)$ eine endliche Untergruppe, die auf $\C^2$ durch Matrixmultiplikation von rechts operiert und damit auf $\C[z_1, z_2]$ durch $(gf)(z_1, z_2)=f((z_1, z_2)g))$ für $f\in \C[z_1, z_2]$ und $g\in G$. Der Invariantenring $\C[z_1, z_2]^G$ ist eine endlich erzeugte $\C$--Algebra, d.h. es gibt endlich viele invariante Polynome $X_1, \ldots, X_n\in \C[z_1, z_2]^G$ mit $X_i(0)=0$, und endlich viele Relationen $f_j(X_1, \ldots, X_n)=0$ mit $f_j\in \C[x_1, \ldots, x_n], j=1, \ldots, k$, so dass die kanonische Abbildung 
\[
\C[x_1, \ldots, x_n]/\langle f_1, \ldots, f_k\rangle\to\C[z_1, z_2]^G, x_i\mapsto X_i,
\]
ein Isomorphismus ist. Der Orbitraum von $G$ wird durch die Bijektion
\[
\C^2/G\to X:=\{x\in \C^n|f_1(x)=\cdots=f_k(x)=0\}
\]
auf kanonische Weise eine normale analytische Menge in $\C^n$. Der Raumkeim $(\C^2/G,0)=(X,0)$ heißt \textit{Quotientensingularität}. Diese Interpretation der $ADE$--Singularität als Quotientensingularität hatte Klein natürlich noch nicht, sie geht hauptsächlich auf Du Val zurück (\cite{DV1957}).

Die Tatsache, dass es sich um Quotientensingularitäten handelte, und natürlich die expliziten Gleichungen waren wesentlich für Brieskorns Beweis der simultanen Auflösung der ADE-Singularitäten. Die Durchführung dieses Beweises war aber nicht geradlinig, sondern wurde unterbrochen von anderen großartigen Entdeckungen Brieskorns. Insbesondere der Fall der $E_8$--Singularität bereitete größere Schwierigkeiten, was auch darin zum Ausdruck kommt, dass Brieskorn den Beweis für die $A_k, D_k, E_6$ und $E_7$ Singularitäten 1966  in \cite{EB1966a} und für $E_8$ erst 1968 in \cite{EB1968a} veröffentlichte.

Um für die simultane Auflösung den Basiswechsel zu bestimmen, muss man die lokale Monodromie analysieren. Da die ADE--Singularitäten gewichtet homogene (oder quasihomogene) Gleichungen haben, ist die geometrische Monodromie analytisch berechenbar, und sie ist von endlicher Ordnung. Im Fall der Quotientensingularitäten vom Typ $A_k, D_k, E_k$ ist die Monodromieoperation auf der mittleren Homologie der Faser ein Coxeter--Element der Spiegelungsgruppe, also das Produkt der zu einer Kammer gehörenden Erzeuger von $G$. Es liegt also nahe, einen Basiswechsel $s=t^d$ zu betrachten, wobei $d$ die Ordnung des Coxeter--Elements, die Coxeter--Zahl, ist. Das Faserprodukt von $X\to S$ und dem Basiswechsel $T\to S$ hat dann die Gleichung
\[
f(x,y,z)-t^d=0,
\]
wobei $f(x,y,z)=0$ die Gleichung einer ADE--Singularität ist. Im $A_1$--Fall von Atiyah ist dies die Gleichung 
\[
x^2+y^2+z^2-t^2=0,
\]
die sich nach Koordinatenwechsel als
\[
z_1z_2-z_3z_4=0
\]
schreiben lässt. Dies ist die Gleichung einer 3--dimensionalen singulären Quadrik $Q_3$ im $\C^4$, also des Kegels über einer nicht--singulären Quadrik im $\P^3$. Durch Aufblasen der Kegelspitze erhält man eine nicht--singuläre Varietät $Y$ über $T$, deren exzeptioneller Divisor isomorph zu $\P^1\times \P^1$ ist und der auf zwei Weisen zu $\P^1$ niedergeblasen werden kann. Die so entstehenden Varietäten $Y_1$ und $Y_2$ sind zwei verschiedene simultane Auflösungen der gegebenen Familie. Sie sind sogenannte  \textit{kleine Auflösungen} der 3-dimensionalen Quadrik $Q_3$, da die exzeptionelle Menge eine (rationale) Kurve ist. Die natürliche bijektive Äquivalenz zwischen $Y_1$ und $Y_2$ heißt (Atiyah--) Flop. Flops und Flips spielen bei den sogenannte minimalen Modellen einer algebraischen Varietät eine fundamentale Rolle, Flips bei der Konstruktion selbst, die bisher nur bis zur Dimension 3 bekannt ist, während verschiedene minimale Modelle durch eine Folge von Flops verbunden sind.

Da für $A_k, E_6$ und $E_8$ die Gleichung $f(x,y,z)-t^d=0$ die Form
\[
x^a+y^b+z^c+t^d=0  
\]
hat, versuchte Brieskorn für diese Singularitäten kleine Modifikationen zu konstruieren, indem er sie auf andere Varietäten abzubilden versuchte, für die eine kleine Auflösung bereits bekannt war. Für den quadratischen Kegel $Q_3$ bedeutete dies, die Gleichung $f(x,y,z)-t^d$ in der Form $\phi_1\phi_2-\phi_3\phi_4$ zu schreiben. Mit solchen Methoden gelang es Brieskorn 1964, simultane Auflösungen der $A_k, D_k, E_6$ und $E_7$--Singularitäen zu konstruieren (veröffentlich in \cite{EB1966a}). Es stellte sich dabei auch heraus, dass eine simultane Auflösung der Abbildung einer 3--Mannigfaltigkeit auf eine 1--Mannigfaltigkeit nur für die ADE--Singularitäten möglich ist. Damit blieb nur noch der Fall der $E_8$--Singularität übrig.

Die $E_8$--Singularität erwies sich jedoch als äußerst hartnäckig und Brieskorn gelang es zunächst nicht, für diese eine simultane Auflösung zu konstruieren. Die verschiedenen Versuche dazu führten ihn jedoch zu überraschenden Ergebnissen über die Topologie und Differentialtopologie von Singularitäten, auf die ich im nächsten Abschnitt eingehe.

Die simultane Auflösung der \textit{Ikosaeder--Singularität} $E_8$ fand Brieskorn im September 1966 (veröffentlicht in \cite{EB1968a}) unter Verwendung sehr klassischer algebraischer Geometrie, wie er selbst schreibt. So verwendete er eine Arbeit von Max Noether von 1889 über rationale Doppelebenen und Eigenschaften exzeptioneller Kurven auf rationalen Flächen, die durch das Aufblasen von 8 Punkten auf einer ebenen Kubik entstehen. Brieskorn fand, dass es rund 700 Millionen simultane Auflösungen von $E_8$ gibt, und zwar genau $2^{14}\cdot 3^5\cdot 5^2\cdot 7$, die Ordnung der Weyl-Gruppe vom Typ $E_8$. Die Divisorenklassengruppe des lokalen Rings der Singularität
\[
x^2+y^3+z^5+t^{30}=0
\]
hat die Struktur des Gitters der Gewichte des Wurzelsystems von $E_8$. Brieskorn konstruierte die kleinen Auflösungen dieser Singularität durch Kurven mit $E_8$ als dualem Graph und damit die simultanen Auflösungen der Flächensingularität $E_8$. Die verschiedenen simultanen Auflösungen entsprechen den Weylkammern, wobei das Aufblasen irgendeiner Idealklasse in jeder Kammer zu einer simultanen Auflösung führt.

Die Untersuchungen zur simultanen Auflösung der ADE--Singularitäten als Quotientensingularitäten von $\C^2$ nach einer endlichen Untergruppe von $\SL(2, \C)$ führten Brieskorn in \cite{EB1968b} zur Untersuchung allgemeiner Quotientensingularitäten $\C^2/G$, wobei $G$ eine beliebige endliche Untergruppe von $\GL(2, \C)$ ist. Er klassifizierte diese Singularitäten unter Verwendung von Ergebnissen von Mumford, Hirzebruch und vor allem von Prill durch Auflistung sämtlicher kleiner (d.h. kein Element von $G$ hat $1$ als Eigenwert mit Multiplizität $1$) Untergruppen $G\subset\GL(2, \C)$ und bestimmte den Auflösungsgraphen von $\C^2/G$, bewertet durch die Schnittmultiplizitäten der exzeptionellen Kurven. Brieskorn zeigte, dass dieser Auflösungsgraph die Singularität bis auf analytische Isomorphie bestimmt, und daraus folgerte er das bemerkenswerte Resultat über die Einzigartigkeit der $2$--dimensionalen Ikosaeder--Singularität:

\textit{Der Ring $\C\{x,y,z\}/\langle x^2+y^3+z^5\rangle$ (und dessen Komplettierung) ist der einzige nicht--reguläre faktorielle $2$--dimensionale analytische lokale Ring.}

In Dimensionen 3 gibt es sowohl unendlich viele faktorielle wie nicht--faktorielle lokale Ringe isolierter Hyperflächensingularitäten (ab Dimension 4 hat man immer Faktorialität).

Brieskorns Arbeiten zur simultanen Auflösung und über die Quotientensingularitäten spielten bei der weiteren Entwicklung der Deformationstheorie rationaler Flächensingularitäten eine wichtige Rolle. Ich erwähne hier nur Oswald Riemenschneider, Jonathan Wahl und in Zusammenhang mit dem Programm der \glqq minimalen Modelle\grqq\ Shigefumi Mori, János Kollár, Miles Reid und Vyacheslav Shokurov.        
\medskip

\textbf{Topologie von Singularitäten und exotische Sphären}

Im September 1965 trat Brieskorn eine C.L.E. Moore Instructorship am MIT in Cambridge/Massachusetts an. Das Problem der simultanen Auflösung von $E_8$ war  zu der Zeit noch nicht gelöst und Brieskorn suchte nach Lösungen in Diskussionen u.a. mit Heisuke Hironaka auf der Arbeitstagung 1965 in Bonn und mit Michael Artin und David Mumford am MIT. Brieskorn versuchte, die Divisorenklassengruppe von 
\[
x^2+y^3+z^5+t^{30}=0
\]
zu berechnen, aber Mumford riet, zunächst die einfachere Gleichung $x^2+y^3+z^5+t^2=0$, also die Gleichung der $3$--dimensionalen $E_8$--Singularität zu untersuchen. Zur Übung startete Brieskorn mit der $3$--dimensionalen $A_2$--Singularität
\[
z_0^3+z_1^2+z_2^2+z_3^2=0
\]
und stellte fest, dass sie faktoriell ist. Danach untersuchte er wieder die $3$--dimensionale $E_8$--Singularität und zeigte durch ziemlich mühsame explizite Auflösung, dass sie ebenfalls faktoriell ist und dass die zweite Kohomologiegruppe des Umgebungsrandes der Singularität verschwindet. Da Brieskorn eigentlich eine nicht--triviale Divisorenklassengruppe erwartet hatte, wandte er sich wieder der $A_2$ zu, um deren Topologie besser zu verstehen. Im September 1965 machte er dann die unerwartete Entdeckung, \textit{dass der Umgebungsrand der $3$--dimensionalen $A_2$--Singularität homöomorph zur $5$--dimensionalen Sphäre ist.} Der Umgebungsrand einer Hyperflächensingularität im $\C^{n+1}$ ist der Durchschnitt mit einer hinreichend kleinen reellen Sphäre im $\R^{2n+2}=\C^{n+1}$ um den singulären Punkt. Für eine isolierte Singularität ist dies eine kompakte $(2n-1)$--dimensionale reell analytische Mannigfaltigkeit. Die Singularität selbst, also die Nullstellenmenge der definierenden Gleichung, ist nach Milnor topologisch der Kegel über dem Umgebungsrand mit dem singulären Punkt als Kegelspitze. Daraus folgt, dass die 3-dimensionale $A_2$--Singularität topologisch eine Mannigfaltigkeit ist. Diese Entdeckung kam völlig überraschend, denn in \cite{DM1961} hatte David Mumford gezeigt, dass isolierte Singularitäten algebraischer Flächen nie topologisch trivial sein können, es sei denn die Singularität ist analytisch nicht--singulär. Brieskorn veröffentlichte dann in \cite{EB1966b}, \textit{dass für alle ungeraden $k\geq 3$ die Singularitäten }
\[
z_0^3+z_1^2+\cdots +z_k^2=0,
\]
\textit{topologische Mannigfaltigkeiten sind}, Mumfords Satz also ein spezielles Phänomen in Dimension zwei ist.

Die Entwicklungen im Jahr 1965/66, die dann zur Entdeckung der exotischen Sphären als Umgebungsränder von Singularitäten führten, sind im Rückblick auch heute noch faszinierend, vor allem wegen des Zusammenspiels der Ideen mehrerer Beteiligter, die durch glückliche Umstände zustande kam. Hirzebruch hat über diese Entdeckung im Seminar Bourbaki \cite{FH1967} und später auf der Singularitätentagung in Oberwolfach aus Anlass von Brieskorns 60. Geburtstag berichtet, eine Kurzfassung findet sich in \cite{FH1996}. 

Hirzebruch berichtete in Rom auf einer Konferenz über Brieskorns simultane Auflösung der Singularitäen vom Type $A_k, D_k, E_6$ und $E_7$, als ihn dort am 28.09.1965 ein Brief von Brieskorn erreichte, in dem dieser schrieb: 

\textit{\glqq Ich habe in den letzten Tagen die etwas verwirrende Entdeckung gemacht, dass es vielleicht $3$--dimensionale normale Singularitäten gibt, die topologisch trivial sind. Ich habe dieses Beispiel heute Nachmittag mit Mumford diskutiert, und er hat bis heute Abend noch keinen Fehler gefunden; hier ist es: $X=\{x\in\C^4|x_1^2+x_2^2+x_3^2+x_4^3=0\}$\grqq.}

Dieses Resultat von Brieskorn war damals ziemlich aufregend und stimulierte Hirzebruch, Milnor und andere zu weiteren Untersuchungen der Topologie isolierter Singularitäten. Natürlich gab es in dieser Zeit keine E-Mail, dafür aber eine umfangreiche Korrespondenz zwischen Brieskorn, Hirzebruch, Jänich, Milnor und Nash. Hirzebruch schrieb im März 1966 an Brieskorn, dass er einen engen Zusammenhang zwischen den Arbeiten von Klaus Jänich, der ebenfalls ein Schüler von Hirzebruch war, über die Klassifikation spezieller $O(n)$--Mannigfaltigkeiten und den von Brieskorn untersuchten Umgebungsrändern von Singularitäten gefunden hätte. Jänich hatte die Operation einer kompakten Liegruppe $G$ auf einer differenzierbaren Mannigfaltigkeit $X$ ohne Rand untersucht. Für spezielle Operationen ist der Orbitraum $X/G$ auf kanonische Weise eine differenzierbare Mannigfaltigkeit mit Rand. Durch Brieskorns Arbeiten motiviert, betrachtete Hirzebruch den Umgebungsrand $\Sigma=\Sigma(k+1, 2, \ldots, 2)$ der $A_k$--Singularitäten im $\C^{n+1}$, der durch die folgenden Gleichungen gegeben ist:
\[
\begin{array}{lcl}
z_0^{k+1}+z_1^2+\cdots +z_n^2 & = & 0\ ,\\
|z_0|^2+|z_1|^2+\cdots +|z_n|^2 & = & 1\ .
\end{array}
\]
Er bewies, dass die orthogonale Gruppe $O(n)$ auf $\Sigma$ auf spezielle Weise im Sinne von Jänich operiert, mit Orbitraum der $2$--dimensionale Scheibe $D^2$, und dass $\Sigma$ für gerades $k$ eine Homologiesphäre ist. Noch aufregender war aber seine Entdeckung, dass $\Sigma$ für $n=5$ und $k=2$ eine exotische $9$--Sphäre ist, also homöomorph aber nicht diffeomorph zur Standardsphäre $S^9$. Der Umgebungsrand $\Sigma(3,2,2,2,2,2)$ der $5$--dimensionalen $A_2$--Singularität stellte sich als die von Michel Kervaire durch sogenanntes \glqq plumbing\grqq\ von zwei Kopien des tangentialen Scheibenbündels der $5$--Sphäre konstruierte $9$--dimensionale exotische Kervaire--Sphäre heraus.

Der Brief von Hirzebruch an Brieskorn vom 24.03.1966, in dem er seine Entdeckung beschreibt, beantwortet Brieskorn am 29.03.1966 u.a. mit folgenden Worten: 

\textit{\glqq Klaus Jänich und ich hatten von diesem Zusammenhang unserer Arbeiten nichts gemerkt, und ich war vor Freude ganz außer mir, wie Sie nun die Dinge zusammengebracht haben.\grqq}

Während in Kervaires Konstruktion die exotische Sphäre eine parallelisierbare Mannigfaltigkeit berandet, ist $\Sigma$ der Umgebungsrand einer Singularität und es blieb zunächst mysteriös, wo die parallelisierbare Mannigfaltigkeit im Singularitätenbild zu finden war. Exotische Sphären waren zuerst von John Milnor in \cite{JM1956} entdeckt worden und diejenigen einer festen Dimension bilden eine abelsche Gruppe $\Theta_n$, mit der zusammenhängenden Summe als Gruppenoperation. Eine wichtige Untergruppe ist $bP_{n+1}$, die aus denjenigen Sphären besteht, die eine parallelisierbare Mannigfaltigkeit beranden. Dies war 1963 von Kervaire und Milnor gezeigt worden. Sie bewiesen außerdem, dass die Gruppe $\Theta_n$ für $n\geq 5$ endlich ist und dass $b P_{4k+2}$ entweder $0$ oder $\Z/2\Z$ ist und dass der zweite Fall genau dann auftritt, wenn der Erzeuger von $bP_{4k+2}$, die $(4k+1)$--dimensionale exotische Kervaire--Sphäre ist. 

Milnor war durch Brieskorns Beispiel des Umgebungsrandes $\Sigma(3,2,2,2)$ als topologische Mannigfaltigkeit angeregt worden, die Umgebungsränder weiterer Singularitäten zu untersuchen und erläuterte seine Überlegungen in einem Brief an John Nash vom April 1966. Brieskorn zitiert in \cite{EB2010} aus diesem Brief:

\textit{\glqq Dear John,}

\textit{I enjoyed talking to you last week. The Brieskorn example is fascinating. After staring at it a while I think I know which manifolds of this type are spheres, but the statement is complicated and the proof doesn't exist yet. Let $\Sigma(p_1, \ldots, p_n)$ be the locus}
\[
z_1^{p_1}+\cdots + z_n^{p_m}=0\ ,\ |z_1|^2+\cdots +|z_n|^2=1\ .\grqq
\]

Dann fährt Milnor mit einer konkreten Vermutung fort, welche dieser Mannigfaltigkeiten topologische Sphären sind. Brieskorn erwähnt weiter, dass der Brief am Rand eine kleine, etwa 1 cm große Skizze enthält, die er aber damals nicht verstanden hätte.

Milnors Skizze, die ich aus \cite{EB2010} kopiert habe (s. Fig. 1) zeigt das Bild der Milnorfaserung und damit die gesuchte parallelisierbare Mannigfaltigkeit. Diese Skizze wurde später in der Singularitätentheorie eine Ikone und zierte fast jeden Vortrag zur Topologie von Singularitäten.
\begin{center}
\includegraphics{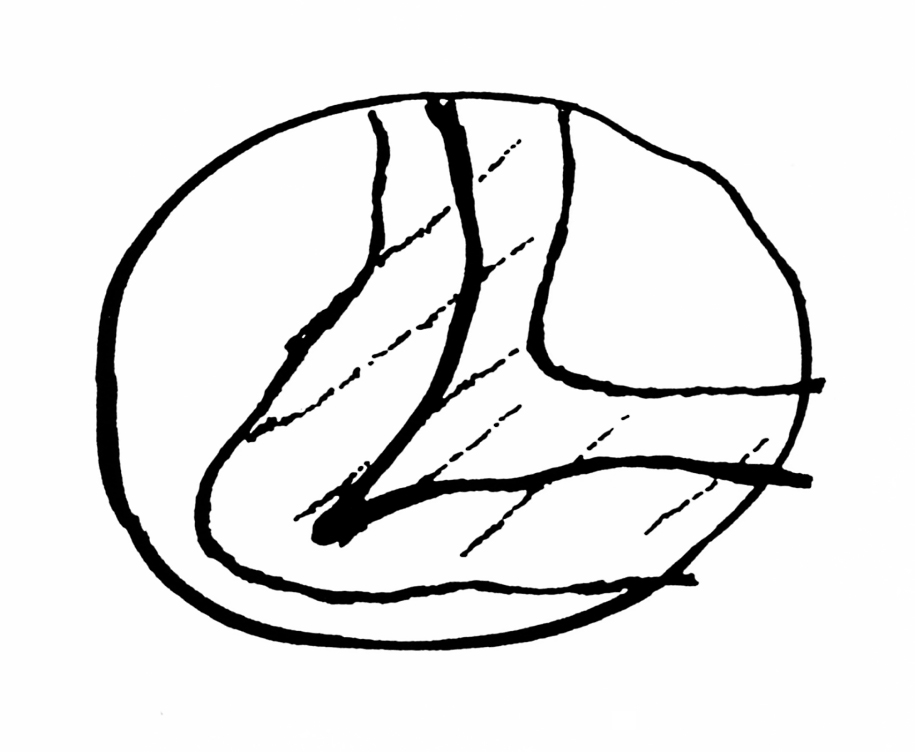}\\
Fig. 1
\end{center}
Man wähle eine hinreichend kleine Kugel $B_\varepsilon$ vom Radius $\varepsilon$ um den isolierten singulären Punkt der Hyperfläche $f(z_1, \ldots, z_n)=0$ und dann eine kleine Kreisscheibe $D_\delta$ vom Radius $\delta\ (<<\varepsilon)$ in der komplexen Ebene um $0$, und Milnors Skizze zeigt $X=f^{-1}(D_\delta)\cap B_{\varepsilon}$. Dann ist $X\smallsetminus f^{-1}(0)$ ein differenzierbares Faserbündel, dessen Faser $X_s=f^{-1}(s), s\neq 0$, als Milnorfaser bezeichnet wird. $X_s$ ist eine $(n-2)$--zusammenhängende parallelisierbare Mannigfaltigkeit, deren Rand diffeomorph zum Rand $\Sigma$ von $X_0$ ist. Also ist die Milnorfaser die parallelisierbare Mannigfaltigkeit, die von der exotischen Kervaire--Sphäre berandet wird und nach der Hirzebruch und Brieskorn gesucht hatten.

Brieskorn gelang es dann in \cite{EB1966c}, die Vermutung von Milnor innerhalb von 14 Tagen vollständig zu beweisen. Gleichzeitig zeigte er durch explizite Berechnung, \textit{dass der Umgebungsrand $\Sigma(2,2,2,3,5)$ der Ikosaeder--Singularität Milnors exotische $7$--Sphäre ist, das Erzeugende der Gruppe $bP_8=\Theta_7$ der Ordnung 28. Alle verschiedenen 28 exotischen differenzierbaren Strukturen auf $S^7$ werden durch den Umgebungsrand $\Sigma(2,2,2,3,6k-1), k=1, \ldots, 28$ gegeben, also durch einfache reell analytische Gleichungen. Darüber hinaus zeigte er, dass jede ungerade dimensionale Sphäre, die eine parallelisierbare Mannigfaltigkeit berandet, diffeomorph zu einem Umgebungsrand $\Sigma(a_1, \ldots, a_m)$ ist.}

Dies galt als Sensation, denn während Milnors erste Konstruktion einer exotischen Sphäre  in der Tat sehr speziell war, war Brieskorns Konstruktion ganz natürlich und alles andere als \glqq exotisch\grqq. 

Dass Brieskorn Milnors Vermutung so schnell beweisen konnte, ist ebenfalls einem glücklichen Umstand zu verdanken. Bei der Durchsicht der neu erschienenen Journale in der Bibliothek des MIT fiel ihm die Arbeit \cite{FP1965} von Frédéric Pham in die Hände, in der dieser, motiviert durch Singularitäten von Feynman--Integralen in der theoretischen Physik, genau die Singularitäten
\[
X_1^{a_1}+\cdots + X_n^{a_n}=0
\]
untersuchte, die auch Milnor in seinem Brief an Nash betrachtet hatte. Pham berechnete für diese Singularitäten den Homotopietyp der Milnorfaser und die Monodromie der Milnorfaserung. Brieskorn nutzte diese Ergebnisse und Hirzebruchs Berechnung der Signatur der Milnorfaser, um die oben erwähnten Resultate zu beweisen. Seitdem werden diese Singularitäten auch {\em Brieskorn-Singularitäten} oder {\em Brieskorn-Pham-Singularitäten} genannt.

Brieskorns Ergebnisse über die exotischen differenzierbaren Strukturen auf den Umgebungsrändern haben zu vielen Anwendungen in Arbeiten anderer Mathematiker zur Differentialtopologie von Mannigfaltigkeiten geführt. In diesem Zusammenhang gibt es aber nur eine Arbeit von Brieskorn selbst, die Konstruktion exotischer Hopf-Mannigfaltigkeiten, zusammen mit Antonius van de Ven in \cite{BV1968}.

Brieskorn nennt die zwei Jahre in Boston und Cambridge als mit die zwei besten seines mathematischen Lebens.
\medskip

\textbf{Picard--Lefschetz--Monodromie und Gauß--Manin--Zusammenhang}

Für eine isolierte Singularität, gegeben durch $f\in\C\{x_0, \ldots, x_n\}, f(0)=0$, betrachten wir Milnors Konstruktion $f\!\!:\!\! X=f^{-1}(S)\cap B\to S$ mit $B=B_\varepsilon$ und $S=D_\delta$ aus dem vorigen Abschnitt. Die nicht--singuläre Milnorfaser $X_s=f^{-1}(s)$ ist eine Deformation von $X_0=f^{-1}(0)$ und zwar die einfachste, da sie durch $f$ selbst gegeben ist. Die Milnorfaser ist die singularitätentheoretische Erklärung dafür, dass die Brieskornschen exotischen Sphären parallelisierbare Mannigfaltigkeiten beranden. Sie erklärt aber noch nicht, wie diese parallelisierbaren Mannigfaltigkeiten durch plumbing konstruiert werden können. Dazu braucht man eine etwas kompliziertere Deformation, eine sogenannte Morsifizierung. Die Idee geht zurück auf die beiden 1897 bzw. 1906 erschienen Bände \glqq Théorie des fonctions algébriques de deux variables indépendantes\grqq\ von Picard-Simart und auf die Monographie \glqq L'analysis situs et la géométrie algébrique\grqq\ von Lefschetz aus dem Jahr 1924. Später wurde daraus die lokale Picard--Lefschetz--Theorie entwickelt, zu der Brieskorn 1970 im Anhang zu \cite{EB1970a} beigetragen hat.

Dort betrachtet Brieskorn eine Deformation
\[
f_a(x)=f(x)-\sum\limits_{i=0}^n a_ix_i\ ,
\]
wobei $a=(a_0, \ldots, a_n)$ hinreichend allgemein gewählt ist. Bezeichnet $\mu=\mu(f)$ die \textit{Milnorzahl} von $f$, also die Vektorraumdimension der Milnoralgebra
\[
\C\{x_0, \ldots, x_n\}/\langle \frac{\partial f}{\partial x_0}, \ldots, \frac{\partial f}{\partial x_n}\rangle,
\]
dann gibt es in der Umgebung von $0\in \C^{n+1}$ genau $\mu$ Punkte $z_r$, so dass $f_a$ in $z_r$ einen gewöhnlichen Doppelpunkt, d.h. eine Singularität vom Typ $A_1$ hat. Eine Deformation mit $\mu$ gewöhnlichen Doppelpunkten nahe $0$ nennt man eine \textit{Morsifikation} von $f$, eine Idee, die auf René Thom zurück geht. Die Abbildung $f_a:X^a=f_a^{-1}(S)\cap B\to S$ ist genau in den Punkten $z_1, \ldots, z_\mu$ singulär und außerhalb der Fasern durch diese Punkte ein differenzierbares Faserbündel mit Faser $X_s^a=f_a^{-1}(s) ,\ s\neq f_a(z_r),$ diffeomorph zur Milnorfaser $X_s$.

Brieskorn betrachtet nun die Milnorsche Konstruktion für die gewöhnlichen Doppelpunkte $z_r$ von $f_a$, also $f_a^r: f_a^{-1}(D_r)\cap B_r\to D_r$, wobei $B_r\subset B$ eine hinreichend kleine Kugel um $z_r$ vom Radius $\rho$ und $D_r\subset S$ eine kleine Scheibe vom Radius $\delta<<\rho$ um $f_a(z_r)$ ist. Dann gilt für geeignete Koordinaten $y_0, \ldots, y_n$ in der Umgebung von $z_r$
\[
f_a^r(y_0, \ldots, y_n)=f_a(z_r)+y_0^2+\cdots y_n^2\ .
\]
Daraus folgt, dass die Milnorfaser $F_a^r$ von $f_a^r$ die $n$--dimensionale Sphäre
\[
S_r^n=\{y|y\text{ reell }, y_0^2+\cdots+y_n^2=\rho\}
\]
als Deformationsretrakt hat. $S_r^n\subset F_a^r$ definiert eine Homologieklasse $d_r$  in $H_n(F_a^r, \Z)$, $r=1, \ldots, \mu$, und dies sind die schon von Lefschetz betrachteten \textit{verschwindenden Zyklen}, da sie sich auf den singulären Punkt $z_r$ zusammenziehen, wenn $\rho$ gegen $0$ geht. Durch Wahl geeigneter Wege $\gamma_r$ in $D$ von einem Randpunkt von $D_r$ zu dem nicht--kritischen Wert $s$ kann man $d_r$ in die Milnorfaser $X_s^a$ transportieren und erhält eine Homologieklasse $e_r$ in $H_n(X_s^a, \Z)$. Brieskorn zeigt dann, dass
\[
H_n(X_s^a, \Z)=\Z e_1\oplus\cdots\oplus\Z e_\mu
\]
gilt.

Dieses Vorgehen liefert für die ADE--Singularitäten die gesuchte Plumbing-Konstruktion der Milnorfaser, und zwar wie folgt. Durch Wahl eines Ehresmannschen Zusammenhangs für das differenzierbare Faserbündel $X^a\smallsetminus \underset{r}{\cup} f_a^{-1}(f_a(z_r))\to D\smallsetminus\{f_a(z_r)|r=1, \ldots, \mu\}$ lassen sich die verschwindenden Zykeln $S_n^r$ selbst über $\gamma_r$ in eingebettete $n$--Sphären in die Milnorfaser $X_s^a$ transportieren, die dort ebenfalls verschwindende Zykeln genannt werden. Diese verschwindenden Zykeln haben Tubenumgebungen in der Milnorfaser, die isomorph zu ihrem tangentialen Scheibenbündel sind. Für die ADE--Singularitäten lassen sich die verschwindenden Zykel so wählen, dass sich die Milnorfaser direkt mit der Plumbing-Konstruktion dieser Scheibenbündel als parallelisierbare Mannigfaltigkeit realisieren lässt.

Das wesentliche Ziel Brieskorns in der Arbeit \cite{EB1970a} war allerdings nicht die Konstruktion verschwindender Zykeln mit Hilfe einer Morsifikation, sondern die Berechnung der algebraischen Monodromie einer isolierten Hyperflächensingularität $f\in \C\{x_0, \ldots, x_n\}$ mit $f\!\!:\!\!X\to S$ wie zu Beginn dieses Abschnitts. Die geometrische Monodromie ist ein Diffeomorphismus der Milnorfaser $X_s$ auf sich selbst, gegeben durch Weg-Lifting eines einfach geschlossenen Weges $\gamma$ um $0$ in $S$ mit Anfangs-- und Endpunkt $s$ in dem Faserbündel
\[
X':=X\smallsetminus X_0\to S\smallsetminus\{0\}=:\!S'\ .
\]
Die geometrische Monodromie induziert die ganzzahlige Monodromie $H^n(X_s, \Z)\xrightarrow{\cong} H^n(X_s, \Z)$ auf der mittleren Kohomologiegruppe der Milnorfaser, die lokale Picard--Leftschetz--Monodromie von $f$, deren charakteristisches Polynom $\Delta_f$ nach Milnor die Topologie des Umgebungsrandes $\Sigma$ von $f$ weitgehend bestimmt.

Brieskorn gibt in dieser Arbeit eine algebraische Beschreibung der \textit{komplexen lokalen Picard--Lefschetz--Monodromie}
\[
h_f:H^n(X_s, \C)\xrightarrow{\cong} H^n(X_s, \C)
\]
und leitet daraus einen Algorithmus zur Berechnung des \textit{charakteristischen Polynoms} $\Delta_f$ ab. 

Zur Berechnung der komplexen Monodromie benutzt Brieskorn holomorphe Differentialformen. Zunächst sind die Kohomologiegruppen $H^p(X_s, \C), s\in S'$, die Fasern eines holomorphen Vektorbündels, deren Garbe der holomorphen Schnitte kanonisch isomorph zu
\[
R^nf_\ast\C_{X'}\otimes_{\C_{S'}} \ko_{S'}
\]
ist. Hierbei ist $R^nf_\ast\C_{X'}$ die $n$--te direkte Bildgarbe der konstanten Garbe $\C_{X'}$. Da sich die Kohomologie der Steinschen Mannigfaltigkeit $X_s$ mittels holomorpher Differentialformen berechnen lässt, betrachtet Brieskorn den Komplex der relativen holomorphen Differentialformen von $X$ über $S$,
\[
\Omega^\bullet_{X/S}=\Omega^\bullet_X/df\wedge\Omega^{\bullet-1} X,
\]
mit dem Differential $\Omega^p_{X/S}\to \Omega^ {p+1}_{X/S}$ induziert durch die äußere Ableitung auf dem Komplex $\Omega^\bullet_X$ der holomorphen Differentialformen auf der Mannigfaltigkeit $X$. Man hat nun einen kanonischen Isomorphismus
\[
R^n f_\ast\C_{X'}\otimes _{\C_{S'}}\ko_{S'} \cong H^n(f_\ast\Omega^\bullet_{X'/S'})
\]
und die rechte Seite, die $n$--te Kohomologiegarbe des Bildgarbenkomplexes $f_\ast\Omega^\bullet_{X'/S'}$, lässt sich durch
\[
\kh^n(X/S):=H^n(f_\ast\Omega^\bullet_{X/S})
\]
auf ganz $S$ fortsetzen. Brieskorn zeigt, dass $\kh^n(X/S)$ auf $S$ kohärent ist und dass für den Halm in $0$ gilt
\[
\kh^n(X/S)_0=H^n(\Omega^\bullet_{X/S, 0})=:H,
\]
diese also nur von der Singularität von $f$ in $0$ abhängt.

$\kh^n(X/S)$ hat als $\ko_S$--Garbe den Rang $\mu=\mu(f)=\dim_\C H^n(X_s, \C)$, $s\in S'$, und Brieskorn vermutete, dass sie lokal frei ist, was kurz darauf von Marcos Sebastiani in \cite{MS1970} bewiesen wurde. Auf $H$ definiert Brieskorn den (meromorphen) \textit{lokalen Gauß--Manin--Zusammenhang} durch die Formel
\[
\triangledown_f\omega=\frac{d\omega}{df}\ .
\]
Das bedeutet, dass für einen Repräsentanten $\widetilde{\omega}\in\Omega^n_{X,0}$ von $\omega$ eine Gleichung $d\widetilde{\omega}=df\wedge\psi$ gilt und $\frac{d\omega}{df}$ die Klasse von $\psi$ in $\Omega^n_{X/S,0}/d\Omega^{n-1}_{X/S, 0}$ bezeichnet. Dass dies wohldefiniert ist, folgt aus dem sogenannten De Rham--Lemma, im Prinzip eine Aussage über die Exaktheit des Koszul--Komplexes für die reguläre Folge $\frac{\partial f}{\partial x_0}, \ldots, \frac{\partial f}{\partial x_n}$. Für $k$ mit $f^k\in\langle\frac{\partial f}{\partial x_0}, \ldots, \frac{\partial f}{\partial x_n}\rangle$ gilt $f^k\frac{d\omega}{df}\in H$ und Brieskorn zeigt:

\textit{$\triangledown_f$ ist ein singulärer (meromorpher) gewöhnlicher Differentialoperator erster Ordnung auf $H$, dessen Monodromie (durch analytische Fortsetzung eines Fundamentalsystems von Lösungen längs eines geschlossenen Weges in $S'$) sich kanonisch mit der lokalen Picard--Lefschetz--Monodromie identifiziert.}

Außerdem beweist Brieskorn, \textit{dass $\triangledown_f$ regulär--singulär ist}, d.h. sich durch eine meromorphe Transformation in einen Differentialoperator mit einem Pol 1. Ordnung überführen lässt. Hieraus leitet Brieskorn einen Algorithmus zur Berechnung des charakteristischen Polynoms $\Delta_f$ der Monodromie von $\triangledown_f$ her.

Da $\Delta_f$ ein ganzzahliges Polynom ist, das in gewissem Sinne algebraisch ist (wie Brieskorn zeigt), folgerte er aus der Regularität von $\triangledown_f$ und der Lösung des 7. Hilbertschen Problems durch Gelfand und Schneider (1934), \textit{dass die Eigenwerte der Monodromie Einheitswurzeln $e^{2\pi i\mu_j}$ mit rationalem $\mu_j$ sind}. Die Aussage dieses Satzes wird auch als \textit{Monodromiesatz} bezeichnet und sie war 1970 schon von Pierre Deligne für globale algebraische Morphismen mit anderen Methoden bewiesen worden. Brieskorns Beweis gilt als besonders elegant.

Die Manuscripta--Arbeit über den lokalen Gauß--Manin--Zusammenhang führte zu bedeutenden Entwicklungen, unter anderem bei Brieskorns Schülern Kyoji Saito, John Peter Scherk, Wolfgang Ebeling, Claus Hertling und dem Autor dieser Zeilen. Ich selbst war Student in Göttingen, als Brieskorn die Arbeit fertigstellte und sein erster Diplomand. Ich erhielt die Aufgabe, die Arbeit auf vollständige Durchschnitte zu verallgemeinern. Die Hauptschwierigkeit war die Verallgemeinerung des Lemmas von De Rham. Mit Hilfe cohomologischer Methoden, die ich auf Anregung von Jean--Pierre Serre während seines Besuchs in Göttingen verwendete, gelang der Beweis und war ein Hauptresultat meiner 1971 fertig gestellten Diplomarbeit. Später wurde das \glqq verallgemeinerte De Rham--Lemma\grqq\ von Kyoji Saito sowie Wolfgang Ebeling und Sabir Gusein--Zade weiter verallgemeinert. Der lokale Gauß--Manin--Zusammenhang zusammen mit dem Index--Satz für regulär--singuläre Differentialoperatoren von Malgrange waren auch der Schlüssel zum Beweis einer algebraischen Formel für die Milnorzahl isolierter Singularitäten vollständiger Durchschnitte in \cite{GG1975} (angekündigt in der gemeinsamen Arbeit \cite{BG1975}), die unabhängig mit topologischen Methoden von Lê D$\tilde{\text{u}}$ng Tráng hergeleitet wurde.

Neben dem Modul $H=H^n(\Omega^\bullet_{X/S,0})$ führte Brieskorn die beiden Moduln
\[
\begin{array}{lcl}
H': & = & df\wedge\Omega^n_{X,0}/df\wedge d\Omega^{n-1}_{X,0}\ \text{ und}\\
H'': & = & \Omega^{n+1}_{X,0}/df\wedge d\Omega^{n-1}_{X,0}
\end{array}
\]
ein, die ebenfalls freie $\ko_{S,0}$--Moduln vom Rang $\mu(f)$ sind, wobei der Gauß--Manin--Zusammenhang dann eine Abbildung $\triangledown_f\!\!:\!\! H'\to H'', [df\wedge\omega]\to [d\omega]$, ist. Die Bedeutung dieser ad hoc Definition war zunächst überhaupt nicht klar, sie sollte sich später aber als grundlegend herausstellen. $H'$ und $H''$ werden heute als \textit{Brieskorn--Gitter} bezeichnet und speziell $H''$ spielt eine wichtige Rolle u.a. bei der Untersuchung der gemischten Hodge--Struktur isolierter Singularitäten und in Kyoji Saitos \glqq higher residue pairings\grqq. Wichtige Arbeiten hierzu stammen neben den schon genannten Schülern von Brieskorn von Morihiko Saito, Daniel Barlet, Claude Sabbah, Mathias Schulze und Christian Sevenheck, um nur einige zu nennen.
\medskip

\textbf{Einfache Singularitäten und einfache Liegruppen}

Brieskorns Arbeit über die simultane Auflösung der einfachen Singularitäten führte zu einem seiner wichtigsten Ergebnisse, dem Zusammenhang zwischen den ADE--Singularitäten und einfachen Liegruppen. Er berichtete auf dem Internationalen Mathematikerkongress 1970 in Nizza darüber und veröffentlichte die Ergebnisse in der kurzen Arbeit  \cite{EB1970b}. 

Alexander Grothendieck hatte Brieskorns Arbeiten über die simultane Auflösung gelesen und wurde dadurch zu Vermutungen angeregt, die er Brieskorn mitteilte. Während Brieskorn $1$--parametrige Deformationen der ADE--Singularitäten studiert hatte, die durch das definierende Polynom gegeben waren, regte Grothendieck an, die semi-universelle Deformation dieser Singularitäten anzusehen. Er vermutete, dass diese durch die adjungierte Quotientenabbildung der einfachen Lie--Algebra vom Typ $A$, $D$ oder $E$ bestimmt ist und dass eine simultane Auflösung der semi-universellen Deformation der Singularitäten vom entsprechenden Typ mit Hilfe der Springer--Auflösung der nilpotenten Varietät gegeben ist. Grothendieck selbst hatte simultane Auflösungen von Singularitäten von adjungierten Quotientenabbildungen studiert und war durch Brieskorns Arbeiten auf den vermuteten Zusammenhang gestoßen.

Sei $G$ eine einfache komplexe (algebraische) Liegruppe, also eine komplex algebraische Mannigfaltigkeit mit regulärer Gruppenaktion, die als Gruppe einfach ist. Ist $G$ einfach zusammenhängend, dann ist $G$ durch seine Lie--Algebra $\mathfrak{g}$ bis auf Isomorphie eindeutig bestimmt. Die einfachen Liegruppen entsprechen dabei den einfachen Lie--Algebren und diese sind durch ihr fundamentales Wurzelsystem klassifiziert. Die Wurzelsysteme wiederum werden durch ihr Dynkin--Diagramm (auch Coxeter--Dynkin--Witt--Diagramm) beschrieben und bestimmen $G$ bis auf Isomorphie. Die Klassifikation aller Dynkin--Diagramme, die so aus den einfachen Liegruppen entstehen, liefert vier unendliche Serien $A_k(k\geq 1), B_k (k\geq 2), C_k (k\geq 3), D_k (g\geq 4)$ und die fünf exzeptionellen Fälle $E_6, E_7, E_8, F_4$ und $G_2$. Die \textit{Dynkin--Diagramme} vom Type $A_k, D_k, E_6, E_7, E_8$ zeichnen sich dadurch aus, dass sie homogen sind, ihre Wurzelsysteme also gleich lange Wurzeln haben. Diese Diagramme haben die folgende Gestalt (ADE--Graph):

$A_k: \includegraphics{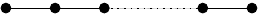}$\hskip1cm $D_k:\includegraphics{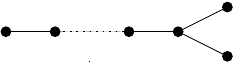}$\ \ \ \ (je $k$ Punkte)\\

$E_6: \includegraphics{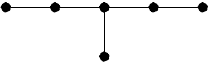}$\hskip1.55cm$E_7: \includegraphics{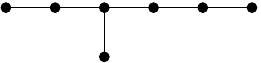}$\hskip 1cm  $E_8: \includegraphics{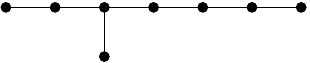}$

Der Name ADE--Singularität für die Quotientensingularitäten nach den endlichen Untergruppen von $\SL(2, \C)$ stammt von der Beziehung zu den einfachen Liegruppen vom Type $A_k, D_k$ oder $E_k$. Auf dem ICM 1970 in Nizza präsentierte Brieskorn die Konstruktion der ADE--Singularitäten sowie ihrer seminuniversellen Deformation direkt mit Hilfe der entsprechenden Liegruppe wie folgt.

Betrachtet man die Operation der einfachen komplexen Liegruppe $G$ auf sich selbst durch Konjugation, so heißt $x\in G$ \textit{regulär}, falls der Orbit von $x$, also seine Konjugationsklasse in $G$, maximale Dimension hat. Ist $d$ diese Dimension, so ist die nächst kleinere Orbitdimension $d-2$ und Elemente dieser Orbitdimension heißen \textit{subregulär}. $G$ besitzt genau einen regulären Orbit, in dessen Abschluss $1\in G$ liegt, und der Abschluss dieses Orbits wird mit Uni $(G)$ bezeichnet, da er die Varietät der unipotenten Elemente der Gruppe ist. Das Komplement des regulären Orbits in Uni $(G)$ hat die Kodimension 2 und ist selbst der Abschluss genau eines subregulären Orbits. Ist $x\in \text{ Uni } (G)$ ein beliebiges Element, dann betrachtet man eine kleine Scheibe $X\subset G$ durch $x$, die transversal zu dem Orbit von $X$ ist, und eine reguläre Projektion $\pi:(G,x)\to (X,x)$ komplexer Raumkeime. Der Raumkeim $(X,x)$ hat in $(G,x)$ komplementäre Dimension zum Orbit von $x$ und ist glatt, falls $x$ reguläres Element ist. Erst bei Schnitten durch nicht reguläre Orbiten entstehen Singularitäten.
Ist $x$ subregulär, so hat $X\cap \text{ Uni } (G)$ die Dimension zwei und eine isolierte Singularität in $x$.

Sei $x=x_sx_n$ die Jordan--Zerlegung von $x\in G$ in einen halbeinfachen und unipotenten Anteil. Ordnet man $x$ die Konjugationsklasse von $x_s$ zu, so erhält man einen Morphismus $\Phi:G\to T/W$, wobei $T$ ein maximaler Torus in $G$ und $W$ die Weyl--Gruppe ist. $\Phi$ heißt die adjungierte Quotientenabbildung. Jede Faser von $\Phi$ ist die Vereinigung endlich vieler Konjugationsklassen, $T/W$ ist eine $k$--dimensionale komplexe Mannigfaltigkeit ($k$ = Anzahl der Ecken des Dynkin--Diagramms) und $\Phi$ bildet Uni$(G)$ auf $1\in T/W$ ab. Mit den eingeführten Bezeichnungen bewies Brieskorn in \cite{EB1970b}. 

Sei $G$ eine einfach zusammenhängende komplexe Liegruppe vom Typ $A_k, D_k, E_6, E_7, E_8$ und $x\in G$ ein subreguläres unipotentes Element. Dann gilt:

\begin{enumerate}
\item [(1)] \textit{$(X\cap \text{ Uni }(G),x)$ ist isomorph zu einer ADE--Singularität vom gleichen Typ wie $G$.}

\item [(2)] \textit{Der adjungierte Quotientenabbildungskeim in $x$ faktorisiert als $\Psi\circ\pi$,
\[
\Phi: (G,x)\xrightarrow{\pi}(X,x)\xrightarrow{\Psi}(T/W, 1),
\]
wobei $\Psi$ die semi-universelle Deformation der entsprechenden Quotientensingularität ist.}
\item [(3)] \textit{Sei 
\[
\xymatrix{\Gamma\ar[r]\ar[d] & G\ar[d]\\
T\ar[r] & T/W}
\]
die von Grothendieck konstruierte simultane Auflösung der adjungierten Quotientenabbildung, mit $\Gamma=\{(x,B)|x\in G, B \text{ Borel--Untergruppe, die } x \text{ enthält}\}$ und $Y$ das Urbild der transversalen Scheibe $X$ in $\Gamma$. Dann folgt:
\[
\xymatrix{Y\ar[r]\ar[d] & X\ar[d]\\
T\ar[r] & T/W
}
\]
ist eine simultane Auflösung der seminuniversellen Deformation der Quotientensingularität vom Typ $A_k, D_k, E_k$.}
\end{enumerate}

Der Beweis von Brieskorn benutzt wesentlich, dass $\Phi$ durch gewichtet--homogene Polynome gegeben ist und dass die ADE-Singularitäten durch ihre Gewichte charakterisiert sind. Die Weylgruppe $W$ ist übrigens gleich der Monodromiegruppe der Singularität. Der von Brieskorn in \cite{EB1970b} skizzierte Beweis wurde von Peter Slodowy vollständig ausgearbeitet. Slodowy erweiterte die Konstruktion später auf alle einfachen Liegruppen, also auch auf die inhomogenen Wurzelsysteme $B_k, C_k, F_4$ und $G_2$, sogar über Körpern beliebiger Charakteristik \cite{PS1980}.

Der von Brieskorn gezeigte Zusammenhang zwischen Liegruppen und Singularitäten führte zu weiteren Untersuchungen. Eine ganz andere Konstruktion der ADE--Singularitäten mit Hilfe der einfachen algebraischen Gruppen vom Typ $A, D$ und $E$ stammt von Friedrich Knop \cite{FK1987}, wobei die Singularitäten allerdings in unterschiedlichen Dimensionen realisiert werden. Eine Klärung des Auftretens der Polyedergruppen in Brieskorns Konstruktion, und damit eine direkte Beziehung zwischen den einfachen Liegruppen und den endlichen Kleinschen Gruppen gelang Peter Kronheimer \cite{PK1990} mit differentialgeometrischen Methoden.  Brieskorn hatte am Ende von \cite{EB1970b} noch geschrieben:

{\em \glqq Thus we see that there is a relation between exotic spheres, the icosahedron and $E_8$. But I still do not see why the regular polyhedra come in.\grqq}

Slodowy untersuchte außerdem kompliziertere Singularitäten und brachte sie mit Kac--Moody--Liealgebren in Verbingung. Kurz vor seinem Tod gelang es ihm, alle einfach elliptischen Singularitäten mit Hilfe der adjungierten Quotientenabbildung der unendlich--dimensionalen Schleifengruppe zu konstruieren. Die Arbeit wurde von Stefan Helmke in \cite{HS2004} vollendet. Slodowy hatte viele Jahre danach gesucht und er berichtete Brieskorn, wenige Tage vor seinem Tod und von seiner Krankheit schon schwer gezeichnet, voller Leidenschaft und Begeisterung davon. Dass die Weiterführung seiner Ideen auch in schwerster Stunde Trost und Freude bereiten konnte, hat Brieskorn zutiefst berührt.
\medskip

\textbf{Verallgemeinerte Zopfgruppen, Milnorgitter und Lorentzsche Raumformen}

Die Konstruktion der semi-universellen Deformation einer ADE-Singularität mit Hilfe der adjungierten Operation der entsprechenden einfachen Liegruppe führte Brieskorn zur Untersuchung von Operationen verallgemeinerter Zopfgruppen und damit zu einer Hinwendung zur Beschäftigung mit diskreten Strukturen isolierter Singularitäten. 

Sei $W$ eine endliche Spiegelungsgruppe, die auf dem reellen endlich dimensionalen Euklidischen Vektorraum $E'$ linear operiert und sei $D'\subset E'$ die Vereinigung der Spiegelungshyperebenen $H'_s$, wobei $s$ ein Element der Menge der Spiegelungen $\Sigma$ in $W$ ist. Brieskorn betrachtet die Komplexifizierung $E$ von $E'$ bzw. $H_s$ von $H'_s$ und $D\subset E$ die Vereinigung der $H_s, s\in\Sigma$. Die Operation von $W$ lässt sich kanonisch auf $E$ erweitern und überführt $E_{reg}=E\smallsetminus D$ in sich. $E_{reg}/W$ ist der Raum der regulären Orbits der endlichen komplexen Spiegelungsgruppe $W$, dessen Fundamentalgruppe Brieskorn in der Arbeit \cite{EB1971} berechnet (s. auch \cite{EB1973}). Er zeigt: 

\textit{Die Fundamentalgruppe $\Pi_1(E_{reg}/ W)$ hat eine Präsentation mit Erzeugern $g_s, s\in\Sigma$, und Relationen}
\[
g_sg_tg_s\cdots=g_tg_sg_t\cdots ,
\] 
\textit{mit $m_{st}$ Faktoren auf beiden Seiten.}

Hierbei ist $(m_{st})$ die Coxeter--Matrix von $W$ mit $m_{st}=$ Ordnung von $st$, und die $g_s$ werden durch eine explizite geometrische Konstruktion angegeben.

Für die symmetrische Gruppe $W\!\!=\!\!S_n\ (=A_{n-1})$  ist die entsprechende Fundamentalgrup\-pe die von Emil Artin, dem Vater von Michael Artin, 1925 eingeführte Zopfgruppe $B_n$, wie von Fox und Neuwirth 1962 und unabhängig von Arnold 1968 bewiesen wurde. Die endlichen irreduziblen Spiegelungsgruppen sind klassifiziert und zerfallen in die Typen $A_k, B_k, D_k, E_6, E_7, E_8, F_4, G_2, H_3, H_4$ und $I_2(m), m=5$ oder $m\geq 7$. Die Fundamentalgruppen der regulären Orbits dieser komplexen Spiegelungsgruppen sind also Verallgemeinerungen der Zopfgruppen.

Der Zusammenhang mit Singularitäten kommt daher, dass für $W$ vom Typ $A_k, D_k, E_6, E_7, E_8$ der Raum $E_{reg}/W$ das Komplement der Diskriminanten im Basisraum der semi-un\-iver\-sellen Deformation der einfachen Singularität vom selben Typ ist. Dies folgt aus Brieskorns Konstruktion in \cite{EB1970b}.

Diese verallgemeinerten Zopfgruppen wurden von Brieskorn und Kyoji Saito zu Ehren von Artin in \cite{BS1972}, \textit{Artin-Gruppen} getauft und unter kombinatorischen Gesichtspunkten untersucht. Unter anderem lösen sie für diese Gruppen das Wort- und das Konjugationsproblem  und bestimmen das Zentrum. Diese Ergebnisse wurden etwa gleichzeitig von Pierre Deligne in \cite{PD1972} erzielt und Deligne beweist darüber hinaus, dass die oben betrachteten Räume $E_{reg}/W$ \textit{Eilenberg--MacLane--Räume} sind, was von Brieskorn in \cite{EB1971} vermutet worden war.

In den im Folgenden erwähnten Arbeiten wendet Brieskorn sich diskreten Invarianten spezieller Klassen von Singularitäten zu. Sei $(X_0, x)\subset(\C^{n+1}, x)$ eine isolierte Hyperflächensingularität und $F\!\!:\!\!X\to S, F(x)=0$, ein geeigneter Repräsentant der semi-univer\-sellen Deformation von $(X_0, x)$. Bezeichne $D$ die Diskriminante von $F$, also die Menge der Punkte $s\in S$ für die die Faser nicht glatt ist, dann ist mit $S'=S\smallsetminus D$ die Einschränkung $F\!:\!X'=F^{-1}(S')\to S'$ eines differenzierbaren Faserbündels mit Faser $X_s$, diffeomorph zur Milnorfaser von $(X_0, x)$. Da $X_s$ den Homotopietyp eines Bouquets von $n$--dimensionalen Sphären hat, ist die mittlere Homologiegruppe $H_n(X_s, \Z)$ frei vom Rang $\mu$, der Milnorzahl von $(X_0, x)$. Ist $n$ gerade, so gibt es auf $H_n(X_s, \Z)$ eine ganzzahlige symmetrische quadratische Form, die Schnittform $< ,>$, und das ganzzahlige Gitter
\[
L=H_n(X_s, \Z)
\]
heißt \textit{Milnorgitter} der Singularität.

Wählt man eine generische komplexe Gerade im affinen $S$ enthaltenden Raum nahe bei $0$, dann ist der Durchschnitt mit $S$ eine kleine komplexe Scheibe $\Delta$, die die Diskriminante $D$ in $\mu$ verschiedenen Punkte $c_1, \ldots, c_\mu$ schneidet. Die Einschränkung von $F$ über $\Delta$ ist eine Morsifikation von $(X_0, x)$, wie sie weiter oben beschrieben wurde. Für $s\in\Delta'=\Delta\smallsetminus\{c_1, \ldots, c_\mu\}$ und eine Wahl von Wegen $\gamma_i$ in $\Delta'$ von $s\in \Delta'$ in die Nähe der $c_i, i=1, \ldots, \mu$, erhält man \textit{verschwindende Zykeln} $e_i\in H_n(X_s, \Z)$ mit $\langle e_i, e_i\rangle=-2$. Die Menge der verschwindenden Zykeln wird mit $\Delta^\ast\subset L$ bezeichnet.

Bei geeigneter Wahl der Wege $\gamma_i$ bilden die $e_1, \ldots, e_\mu$ eine Basis des Gitters $L$, die dann eine \textit{ausgezeichnete Basis} genannt wird. Die Menge aller ausgezeichneten Basen von $L$ wird mit $B^\ast$ bezeichnet. 

Für jede Basis $B\in B^\ast$ beschreibt die Matrix der Skalarprodukte der Basiselemente die Bilinearform auf $L$, die durch einen Graphen $D_B$ beschrieben wird. Die Ecken $\{1, \ldots, \mu\}$ von $D_B$ entsprechen den Basiselementen $e_1, \ldots, e_\mu$ und zwei Ecken $i, j$ werden mit $|\langle e_i, e_j\rangle|$ Kanten verbunden, jeweils mit dem Vorzeichen $\pm 1$ von $\langle e_i, e_j\rangle\in\Z$. $D_B$ heißt \textit{(Coxeter--)Dynkin--Diagramm} von $B$ und die Menge aller Dynkin--Diagramme wird mit $D^\ast$ bezeichnet.

Auf $B^\ast$ und damit auf $D^\ast$ existiert eine natürliche Operation der klassischen Zopfgruppe $B_\mu$ mit $\mu$ Strängen, die sich durch elementare Operationen auf dem Niveau der Wege beschreiben lässt. An verschiedenen Stellen weist Brieskorn darauf hin, dass das Verständnis dieser Operation für ein Verständnis der semi-universellen Deformation essentiell sein sollte.

Eine weitere Invariante ist die \textit{Monodromiegruppe} von $(X_0, x)$. Definitonsgemäß ist dies das Bild unter dem kanonischen Homomorphismus von der Fundamentalgruppe des Komplements der Diskriminante in die Automorphismusgruppe des Gitters $L$. Sie wird bereits von den Automorphismen erzeugt, die zu einer ausgezeichneten Basis gehören.

Eine Übersicht über diese Invarianten und die Beziehungen untereinander gibt Brieskorn in dem Übersichtsartikel \cite{EB1983a} und er betont deren Bedeutung für das Verständnis der Geometrie der semi-universellen Deformation. Wichtige Arbeiten zu diesen Invarianten sind die grundlegenden Arbeiten von Andrei Gabrielov  \cite{AG1974} und von Sabir Gusein-Zade \cite{SG1977} sowie die Lecture Notes von Wolfgang Ebeling \cite{WE1987}.

Ein erster Schritt zum Verständnis besteht darin, die Deformationsrelationen zwischen Singularitäten einer festen Modalitätsklasse zu verstehen. Denn wenn eine Singularität in eine andere deformiert, induziert dies eine Inklusion der entsprechenden Milnorgitter. Die Klassifikation der isolierten Hyperflächensingularitäten bezüglich ihrer Modalität (i.e. die  Anzahl der unabhängigen Parameter (Moduln) von Isomorphieklassen in einer Umgebung des Ursprungs der semi-universellen Deformation) wurde von V.I. Arnold in \cite{VA1972} initiiert und ist einer der Startpunkte der Singularitätentheorie mit weitreichenden Ergebnissen. Die Adjazenzen (Deformationsbeziehungen) zwischen den ADE-Singularitäten wurden schon von Arnold bestimmt. In der Arbeit \cite{EB1979} berechnet Brieskorn alle möglichen Adjazenzen innerhalb Arnolds Liste der unimodularen Singularitäten, das ist nach den ADE-Singularitäten die nächst kompliziertere Klasse in Arnolds Hierarchie. Diese Arbeit wird in \cite{EB1983b} verfeinert, wo Brieskorn eine sehr detaillierte Beschreibung der Milnorgitter der 14 exzeptionellen unimodularen Singularitäten gibt.

Die Deformationsbeziehungen innerhalb der unimodularen Singularitäten wurde von Brieskorn mit einer Theorie in Verbindung gebracht, die scheinbar weit entfernt von der Singularitätentheorie ist, nämlich die Theorie der partiellen Kompaktifizierungen beschränkter symmetrischer Gebiete. Bezeichnet $\kf(L)$ den Isotropie--Fahnen-Komplex von $L$, ein Gebäude im Sinne von Tits, dann operiert die Monodromiegruppe $\Gamma$ auf $\kf(L)$ und der $1$--dimensionale simpliziale Komplex $\kf(L)/\Gamma$ ist endlich. Brieskorn beweist in \cite{EB1981} für die einfachsten exzeptionellen Singularitäten $E_{12}, Z_{11}, Q_{10}$, dass sich die Baily--Borel--Kompaktifizierung von $\kf(L)/\Gamma$ mit dem $\C^\ast$--Quotienten des punktierten negativ--graduierten Teil des Basisraumes der semi-universellen Deformation identifizieren lässt. Unabhängig von Brieskorn hat Eduard Looijenga in \cite{EL1983} diese Ergebnisse für alle Dreiecks--Singularitäten $T_{p, q, r}$ bewiesen und später in allgemeinerem Rahmen wichtige neue Kompaktifizierungen lokal symmetrischer Varietäten konstruiert.

In der Arbeit \cite{EB1988} gibt Brieskorn einen Überblick über die Operation der Zopfgruppe auf der Menge $B^\ast$ der ausgezeichneten Basen einer isolierten Singularität. Außerdem führt er das Konzept einer automorphen Menge ein, das viele Aspekte der Operation der Zopfgruppe vereinheitlicht, und das später von mehreren Autoren aufgegriffen und verallgemeinert wurde. Das folgende Zitat aus der Einleitung zeigt noch einmal die Freude Brieskorns an der Einheit der Mathematik, die im Zusammenspiel vieler verschiedener Gebiete der Mathematik zum Ausdruck kommt.

\textit{\glqq The beauty of braids is that they make ties between so many different parts of mathematics, combinatorial theory, number theory, group theory, algebra, topology, geometry and analysis, and, last but not least, singularities\grqq.}

Damit bin ich fast am Ende des Rückblicks auf Brieskorns mathematisches Werk. Zu erwähnen sind natürlich noch seine beiden Lehrbücher \glqq Lineare Algebra und analytische Geometrie I und II\grqq, die auch wegen der historischen Anmerkungen von Erhard Scholz außerordentlich lesenswert sind, sowie \glqq Ebene algebraische Kurven\grqq\ zusammen mit Horst Knörrer, dessen vorerst letzter Nachdruck in englischen Übersetzung 2012 erschien.

Es gibt noch eine mathematische Arbeit aus dem Jahr 2003 \cite{EB2003} zusammen mit seiner Schülerin Anna Pratoussevitch und seinem Schüler Frank Rothenhäusler, die ich erwähnen möchte. Die Ursprünge dieser Arbeit reichen mindestens bis ins Jahr 1992  zurück, als Brieskorns Student Thomas Fischer ein Polyeder entdeckte, das in gewissem Sinne das klassische Dodekaeder verallgemeinert. Es hat eine sehr ähnliche kombinatorische Struktur wie das Dodekaeder, allerdings mit einer Symmetrieachse der Ordnung 7 statt 5. Sei $\Gamma$ eine diskrete Untergruppe der Liegruppe $\widetilde{SU}(1,1)$, die durch Linkstranslationen operiert. Der Quotient $\Gamma\backslash\widetilde{SU}(1,1)$ ist eine \textit{Lorentzsche Raumform} und wird durch ein Fundamentalgebiet $F$ für $\Gamma$ beschrieben. Für kokompaktes $\Gamma$ ist nach Resultaten von Dolgachev $\Gamma\backslash\widetilde{SU}(1,1)$ der Umgebungsrand einer quasihomogenen Flächensingularität und die quasihomogenen Singularitäten von Arnolds Serien $E_k, Z_k$ und $Q_k$ sind von diesem Typ.

Die Autoren beschreiben in dieser Arbeit die Fundamentalgebiete für die entsprechenden Gruppen $\Gamma$ als Polyeder mit total geodätischen Flächen im 3--dimensionalen Lorentz-Raum wobei jede Serie ein reguläres charakteristisches kombinatorisches Muster zeigt, das mit den klassischen Polyedern in Zusammenhang steht. 

Brieskorn beschreibt in dem Film \glqq Science Lives: Egbert Brieskorn\grqq, siehe \cite{EB2010}, die große Freude, die er empfand, als Fischer das $E_{12}$--Polyeder entdeckt hatte und als Pratoussevitch dies auf die unendlichen Serien $E_k, Z_k$ und $Q_k$ ausdehnen konnte. In der Sequenz \glqq Melencolia\grqq\ des Films erläutert er, dass der richtige Beginn dieser unendlichen Serien das Polyeder auf Dürers berühmter Radierung \textit{Melencolia I} sein sollte. In derselben Sequenz geht er auch darauf ein, wie wichtig ihm insbesondere bei der Unterrichtung der Studierenden die Anschauung ist, im Gegensatz zu einer rein analytischen und strukturellen Vorgehensweise. Fischers Entdeckung erfreute Brieskorn so sehr, dass er selbst eine graphische Darstellung davon berechnete und zeichnete und es sein \glqq Opus 2\grqq\ nannte. Er beauftragte einen Künstler mit der Herstellung einer 3D--Plastik des $E_{12}$--Polyeders aus Messing, die er seinem Lehrer Friedrich Hirzebruch zum 75. Geburtstag schenkte. Soweit ich weiß, ließ er nur zwei oder drei Exemplare anfertigen, wovon er mir ein Exemplar schenkte, was mich außerordentlich gefreut hat. Ein Bild dieses Exemplars, das den Kreis von den Anfängen der platonischen Körper und Quotientensingularitäten bis zu den Lorentzschen Raumformen schließt, ist vielleicht ein geeigneter Abschluss für diese Übersicht über Brieskorns mathematisches Werk.

\begin{center}
\includegraphics[width=6cm,height=6cm]{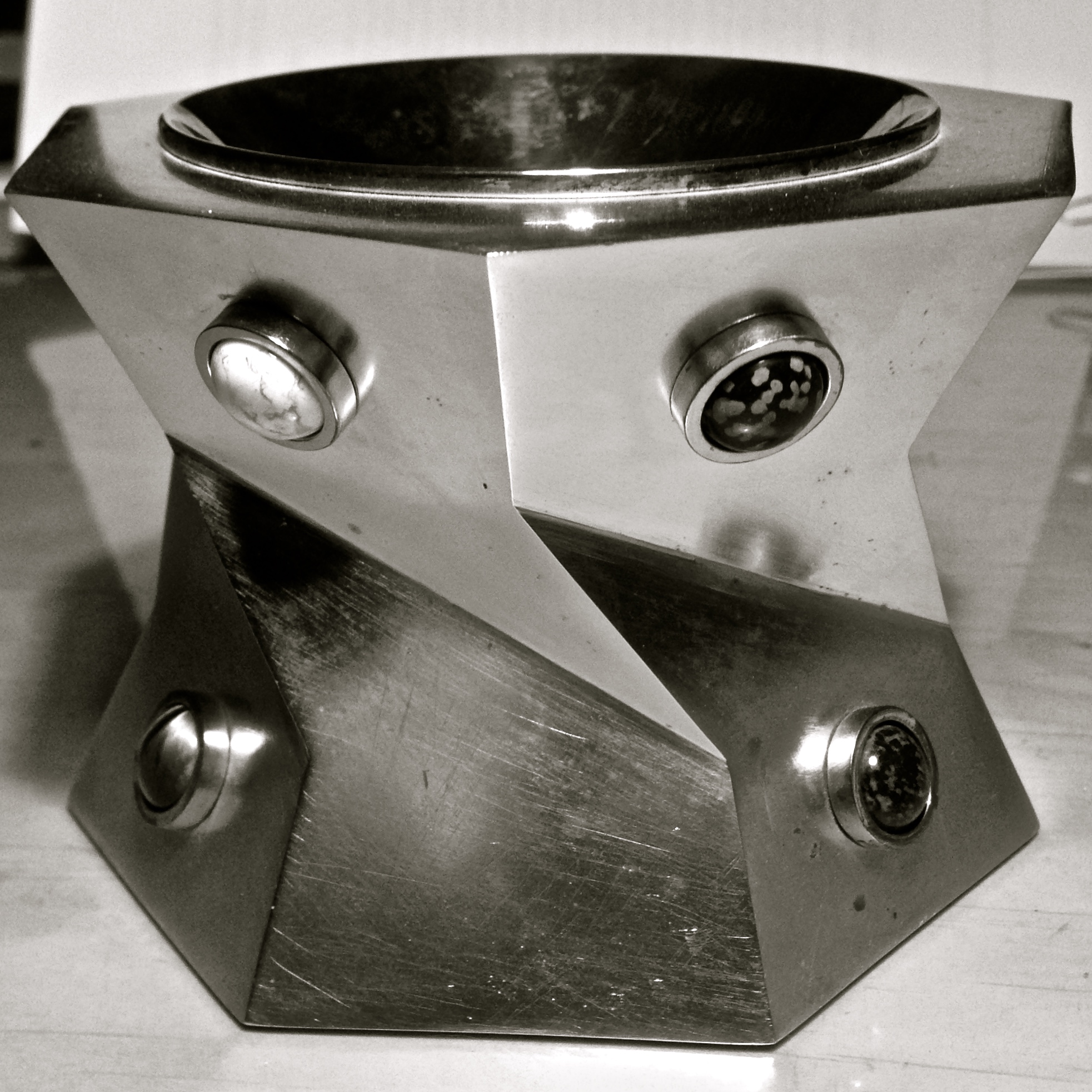}\\
$E_{12}$--Polyeder
\end{center}

      Alle, die Egbert Brieskorn kannten, schätzten sein umfassendes mathematisches Wissen, seine ungeheuer breite allgemeine Bildung, seine scharfe Intelligenz, seine Geradlinigkeit und absolute intellektuelle Redlichkeit, seine freundliche Aufmerksamkeit und Hilfsbereitschaft und seinen besonnenen Rat. Und alle, die ihn näher kennenlernen durften, wissen: Er war ein herzensguter Mensch. Die Geschichte der Wissenschaft wird sein Werk und seinen Namen bewahren, und alle, die ihn kannten und schätzten, werden ihn nicht vergessen.\\

\textbf{\Large Habilitationen}

Brieskorns hat 24 Dissertationen betreut, sie finden sich im "`Mathematics Genealogy Project"'. Von seinen Doktoranden haben sich bisher sieben habilitiert: 

1975:  HAMM, HELMUT AREND: Zur analytischen und algebraischen Beschreibung  der  Picard-Lefschetz-Monodromie. Göttingen.

1980:  GREUEL, GERT-MARTIN: Kohomologische Methoden in der Theorie 
	isolierter Singularitäten. Bonn.

1984:  SLODOWY, PETER: Singularitäten, Kac-Moody-Lie-Algebren.
	assoziierte Gruppen und Verallgemeinerungen. Bonn.

1985:   KNÖRRER, HORST: Geometrische Aspekte integrabler Hamiltonscher 
	Systeme. Bonn.

1986:  EBELING, WOLFGANG: Die Monodromiegruppen der isolierten 	
	Singularitäten vollständiger Durchschnitte.Bonn.

1986:  SCHOLZ, ERHARD: Symmetrie – Gruppe – Dualität. Studien zur
	Beziehung zwischen theoretischer Mathematik und Anwendungen
	in Kristallographie und Baustatik im 19. Jahrhundert. Wuppertal.

2000:  HERTLING, CLAUS: Frobenius-Mannigfaltigkeiten, Gauß-Manin-
	Zusammenhänge und Modulräume von Hyperflächensingularitäten.
	Bonn.

\bigskip

\end{document}